\documentclass[final,onefignum,onetabnum]{siamart190516} 

\usepackage{amssymb}
\usepackage{amsmath}
\usepackage{xparse, mathtools}
\usepackage{tabularx}

\usepackage{hyperref}

\usepackage{caption}
\usepackage{subcaption}

\newcommand{\rev}[1]{#1}
\newcommand{\rone}[1]{#1}
\newcommand{\rtwo}[1]{#1}
\newcommand{\rboth}[1]{#1}
\newcommand{\rme}[1]{#1}

\usepackage{lineno}

\usepackage{bbm}

\newcommand\numberthis{\addtocounter{equation}{1}\tag{\theequation}}

\makeatletter
\DeclareRobustCommand{\pder}[1]{%
  \@ifnextchar\bgroup{\@pder{#1}}{\@pder{}{#1}}}
\newcommand{\@pder}[2]{\frac{\partial#1}{\partial#2}}
\makeatother

\makeatletter
\DeclareRobustCommand{\der}[1]{%
  \@ifnextchar\bgroup{\@der{#1}}{\@der{}{#1}}}
\newcommand{\@der}[2]{\frac{d#1}{d#2}}
\makeatother

\newcommand{\RR}{\mathbb{R}} 

\newcommand{\overbar}[1]{\mkern 1.5mu\overline{\mkern-1.5mu#1\mkern-1.5mu}\mkern 1.5mu}

\newcommand*\diff{\mathop{}\!\mathrm{d}}

\ExplSyntaxOn

\NewDocumentCommand \rvec { s o m }
 {
  \IfBooleanTF {#1}
   { \vectaux*{#3} }
   { \IfValueTF {#2} { \vectaux[#2]{#3} } { \vectaux{#3} } }
  ^\top
 }

\DeclarePairedDelimiterX \vectaux [1] {\lbrack} {\rbrack}
 { \, \dbacc_vect:n { #1 } \, }

\cs_new_protected:Npn \dbacc_vect:n #1
 {
  \seq_set_split:Nnn \l_tmpa_seq { , } { #1 }
  \seq_use:Nn \l_tmpa_seq { \enspace }
 }
\ExplSyntaxOff

\headers{Constrained Pressure-Temperature Residual (CPTR)}{T. Roy, T. B. J\"onsth\"ovel, C. Lemon, A. J. Wathen}

\title{A constrained pressure-temperature residual (CPTR) method for non-isothermal multiphase flow in porous media\thanks{May 27th, 2020. 
\funding{This publication is based on work partially supported by the EPSRC Centre for Doctoral Training in Industrially Focused Mathematical Modelling (EP/L015803/1) in collaboration with Schlumberger.}}}

\author{Thomas Roy\thanks{Mathematical Institute, University of Oxford, Oxford, United Kingdom
  (\email{thomas.roy@maths.ox.ac.uk}, \email{andy.wathen@maths.ox.ac.uk}).}
\and Tom B. J\"onsth\"ovel \thanks{Schlumberger Norway Technology Center, Asker, Norway
  (\email{tjonsthovel@slb.com}).}
\and Christopher Lemon \thanks Schlumberger Abingdon Technology Center, Abingdon, United Kingdom (\email{clemon@slb.com}).
\and Andrew J. Wathen\footnotemark[2]
}

\ifpdf
\hypersetup{
  pdftitle={CPTR},
  pdfauthor={T. Roy, T. B. J\"onsth\"ovel, C. Lemon, A. J. Wathen}
}
\fi

\begin{document}


\maketitle

\begin{abstract}
    For both isothermal and thermal petroleum reservoir simulation, the Constrained Pressure Residual (CPR) method is the industry-standard preconditioner. This method is a two-stage process involving the solution of a restricted pressure system. While initially designed for the isothermal case, CPR is also the standard for thermal cases. However, its treatment of the energy conservation equation does not incorporate heat diffusion, which is often dominant in thermal cases. In this paper, we present an extension of CPR: the Constrained Pressure-Temperature Residual (CPTR) method, where a restricted pressure-temperature system is solved in the first stage. In previous work, we introduced a block preconditioner with an efficient Schur complement approximation for a pressure-temperature system. Here, we extend this method for multiphase flow as the first stage of CPTR. The algorithmic performance of different two-stage preconditioners is evaluated for reservoir simulation test cases.
\end{abstract}

\begin{keywords}
  preconditioning, iterative methods, porous media, thermal reservoir simulation, multiphase flow
\end{keywords}

\begin{AMS}
  	65F08, 65F10, 76S05, 65M08, 65M60
\end{AMS}

\section{Introduction}
\label{sec:intro}
Models of multiphase flow in porous media are used in many applications such as oil and gas recovery, geothermal energy, and carbon sequestration. In some cases, the effects of temperature changes on the fluid flow are significant enough that they are included in the model. Fluid properties affected by temperature include density, viscosity, and phase changes. In the case of heavy oil recovery, non-isothermal flow is an important part of enhanced oil recovery techniques. Under normal reservoir conditions, highly viscous oils do not flow easily to production wells. To reduce the viscosity and thus increase the fluid flow, the temperature inside the reservoir is increased, for instance with the injection of a hot fluid such as steam, or with microwave heating. Decision making regarding the choice of extracting techniques can be guided by thermal reservoir simulation (see for example \cite{coats1980situ}).

In isothermal models, a pressure variable couples a number of secondary variables (saturations/concentrations) that characterize the location of different phases and hydrocarbons. The resulting PDE system is essentially elliptic with respect to pressure and hyperbolic with respect to the secondary variables. Since pressure drives the flow, successful solution techniques usually include a specific treatment of the pressure variable. This is true of classical methods such as the Implicit Pressure Explicit Saturation (IMPES) method, and industry-standard ones such as the Constrained Pressure Residual method (CPR) \cite{wallis1983incomplete, wallis1985constrained}. A preconditioner for the fully implicit system, CPR is a two-stage process where a restricted pressure system is solved, followed by an approximate solution of the coupled system. This method was later improved in \cite{lacroix2003iterative} with the use of Algebraic Multigrid (AMG) \cite{ruge1987algebraic} as a solver for the pressure subsystem. The second stage is usually an incomplete factorization method such as incomplete LU factorization (ILU) \cite{meijerink1977iterative}.

In non-isothermal models, a conservation of energy equation is added to the system along with a temperature (or enthalpy) variable. For fully implicit formulations, the industry-standard preconditioner is CPR where temperature variables are grouped with the secondary variables. This is often appropriate since heat is \rone{advected} similarly to the saturations. However, heat is also diffused through rock and fluids. Diffusion can dominate in cases where the fluid flow is slow, for instance before viscous oils are properly heated, but also due to mesh refinements. In those cases, the second stage of CPR struggles to capture the heat diffusion, and so incomplete factorizations with additional fill are needed. Of course, this remedy is not ideal in terms of scalability and memory requirements.

To avoid having to solve the fully coupled linearized system, there has been some interest in developing alternatives to the standard Newton's method. For example, the Sequential Fully Implicit (SFI) method \cite{jenny2006adaptive} solves flow and transport problems sequentially at the nonlinear level. Although most recent work relates to coupling geomechanics, there is potential for the thermal problem \cite{wong2018sequential}. While the computational cost can be reduced by solving the smaller decoupled systems, considerable effort is needed for such methods to be robust.

An alternative to using AMG only in the first stage of CPR is the use of AMG for the fully coupled system. AMG methods for systems of PDEs \cite{clees2005amg} have been applied successfully to several problems including reservoir simulation \cite{gries2014preconditioning}. In the context of thermal simulation, construction of the grid hierarchy using both pressure and temperature variables has been discussed in \cite{gries2015system}. However, the community's recent focus is again in coupling fluid flow and geomechanics, rather than thermal flow.

In this paper, we present an extension of CPR for non-isothermal flow. Instead of solving a restricted pressure system in the first stage of CPR, we solve a restricted pressure-temperature system, resulting in a Constrained Pressure-Temperature Residual method (CPTR). In \cite{roy2019block}, we presented a Schur complement approximation for the pressure-temperature system for non-isothermal single phase flow. Such a Schur complement approximation leads to an effective block preconditioner. Here, we propose an extension of that method to the multiphase flow situation and use it for the pressure-temperature subsystem in the first stage of CPTR. As an alternative, we also consider applying an unknown-based AMG method to the pressure-temperature subsystem. \rme{In a recent manuscript \cite{cremon2019multi}, promising results were shown for versions of CPTR for thermal-compositional-reactive flow.}

In Section \ref{sec:problem}, we describe a system of PDEs for non-isothermal multiphase flow in porous media as well as its discretization. Then, in Section \ref{sec:precons}, we describe the CPR preconditioner and present our extension, CPTR, as well as an approximate Schur complement solver used in its first stage. Numerical results are presented in Section \ref{sec:results}. We compare the different methods in terms of algorithmic performance, scalability, and robustness to strong cross-couplings. In Section \ref{sec:conclusion}, we conclude with a brief discussion of the results and further work.


\section{Problem statement}
\label{sec:problem}
In this section, we begin by describing a coupled PDE system for thermal multiphase flow, followed by its discretization.
\subsection{Multiphase thermal flow in porous media}
\label{sec:thermaleqns}

We describe the equations for multiphase phase flow in porous media coupled with thermal effects. The number of fluid phases is typically two or three, for example, oil, water, and gas.

The volume fraction occupied by each phase $\alpha$ is called the saturation of that phase, denoted $S_\alpha$. We have the following saturation constraint:
\begin{equation}\label{eq:satconstraint}
 \sum_{\alpha} S_{\alpha} = 1.
\end{equation}
In the absence of phase change, saturations only vary when one phase displaces the other through fluid flow. The ability of a phase $\alpha$ to move is represented by the relative permeability $k_{r\alpha}$, a dimensionless function of $S_\alpha$. 

\subsubsection{Conservation of mass}
For each phase $\alpha$, we have a conservation of mass equation of the following form:
\begin{equation}\label{eq:continuity1}
 \phi\pder{(S_\alpha \rho_\alpha)}{t} + \nabla \cdot (\rho_\alpha \mathbf{u}_\alpha) = f_\alpha \quad \text{in } \RR_+\times\Omega,
\end{equation}
where $\phi$ is the porosity of the rock, $\rho_\alpha$ is the density of the fluid, $\mathbf{u}_\alpha$ is the fluid velocity, $f_\alpha$ is a source/sink term, and $\Omega$ is the spatial domain. \rtwo{We consider a fixed rock strata so that the porosity field $\phi$ is constant in time.} We further assume that the velocity follows Darcy's law \cite{darcy1856fontaines, muskat1936flow}, i.e.
\begin{equation}\label{eq:darcy}
  \mathbf{u}_\alpha = -\mathbf{K}\frac{k_{r\alpha}}{\mu_\alpha}(\nabla p - \rho_\alpha \mathbf{g}),
\end{equation}
where $p$ is the pressure, $\mathbf{K}$ is a permeability tensor, $\mu_\alpha$ is the viscosity, and $\mathbf{g}$ is gravitational acceleration. Note that we have ignored the effect of capillary pressures so that we only have one pressure $p$ for all phases. Then, \eqref{eq:continuity1} becomes
\begin{equation}\label{eq:consmass}
 \phi\pder{(S_\alpha\rho_\alpha)}{t} - \nabla \cdot \left(\rho_\alpha \frac{\mathbf{K}k_{r\alpha}}{\mu_\alpha}(\nabla p - \rho_\alpha \mathbf{g})\right) = f_\alpha \quad \text{in } \RR_+\times\Omega.
\end{equation}
For generality, we first assume mixed and Dirichlet boundary conditions
 \begin{equation}\label{eq:massNeumann}
  -\frac{\mathbf{K}k_{r\alpha}}{\mu_\alpha}(\nabla p - \rho_\alpha \mathbf{g})\cdot \mathbf{n} = g_{N,\alpha} \quad \text{on } \Gamma_N,\quad\text{for all }\alpha,
 \end{equation}
  \begin{equation}\label{eq:massDirichlet}
S_\alpha = g_{D,\alpha}, \;\text{for all }\alpha, \quad p = g_D  \quad\text{on } \Gamma_D, 
 \end{equation}
where $g_{N,\alpha}$ is Neumann boundary data for the flow of phase $\alpha$, $g_D$ and $g_{D,\alpha}$ are Dirichlet boundary data for pressure and saturation $\alpha$, respectively, $\mathbf{n}$ is the unit outward normal vector on $\partial \Omega = \Gamma_N \cup \Gamma_D$, and $\Gamma_D \cap \Gamma_N = \varnothing$. 

\subsubsection{Conservation of energy}

The above equations of fluid flow could be considered on their own as an isothermal model. We are interested in non-isothermal flow, and therefore, in addition to the mass equations, we add a conservation of energy equation for heat, as well as a temperature variable $T$. (Note that enthalpy could also be considered as an independent variable.)

Similarly to the saturations, heat energy is transported by Darcy flow. In addition, however, temperature is also conducted through the porous rock and the different fluids. The conservation of energy is described by the following advection-diffusion equation:
\begin{multline}\label{eq:energy1}
 \pder{}{t} \left(\phi \sum_{\alpha} \rho_\alpha S_\alpha c_{v,\alpha} T + (1-\phi) \rho_r c_r T\right) + \nabla \cdot \sum_{\alpha} \rho_\alpha  c_{v,\alpha}T\mathbf{u}_\alpha - \nabla \cdot (k_T \nabla T) = f_T \\ \text{in } \RR_+\times\Omega,
\end{multline}
where $c_{v,\alpha}$ and $c_r$ are the specific heat of the phase $\alpha$ and the rock, respectively, $\rho_r$ is the density of the rock, $k_T$ is the thermal conductivity field, and $f_T$ is a source/sink term. In this model, the enthalpy of phase $\alpha$ and the rock are given by the simple relations $c_{v,\alpha} T$, and $c_r T$, respectively.

The conductivity depends on the properties of the rock and fluids. We consider the simple relation $k_T = (1-\phi) k_{T,r} + \phi \sum_\alpha S_\alpha k_{T,\alpha}$, where $k_{T,r}$ and $k_{T,\alpha}$ are the conductivities of the rock and the phase $\alpha$, respectively. Assuming Darcy's law, we get
\begin{multline}\label{eq:energy3}
 \pder{}{t} \left(\phi \sum_{\alpha} \rho_\alpha S_\alpha c_{v,\alpha} T + (1-\phi) \rho_r c_r T\right) - \nabla \cdot \sum_{\alpha} \rho_\alpha  c_{v,\alpha}T\mathbf{K}\frac{k_{r\alpha}}{\mu_\alpha}(\nabla p - \rho_\alpha \mathbf{g})\\ - \nabla \cdot (k_T \nabla T) = f_T \quad  \text{in } \RR_+\times\Omega,
\end{multline}
We also first assume Neumann and Dirichlet boundary conditions
\begin{equation}\label{eq:energy4}
  -\left(\sum_{\alpha} \rho_\alpha  c_{v,\alpha}T\mathbf{K}\frac{k_{r\alpha}}{\mu_\alpha}(\nabla p - \rho_\alpha \mathbf{g}) + k_T\nabla T \right) \cdot \mathbf{n} = g^T_N \text{ on } \Gamma^T_N,\; T = g^T_D \text{ on } \Gamma^T_D,
\end{equation}
where $g_N^T$ is Neumann boundary data, $g_D^T$ is Dirichlet boundary data, $\partial \Omega = \Gamma^T_N \cup \Gamma^T_D$, and $\Gamma^T_D \cap \Gamma^T_N = \varnothing$.

Note that $c_{v,\alpha} T$ in the first term of \eqref{eq:energy1} is actually the internal energy of the fluid. For the numerical tests in this paper, we consider liquids and therefore take the internal energy as equal to the enthalpy. For gases, it is more common for internal energy to also have some dependence on pressure and density.


\subsubsection{Coupled problem}
We assume that the relative permeability $k_{r\alpha}$, the density $\rho_\alpha$, and the viscosity $\mu_\alpha$ of phase $\alpha$ are empirically determined functions of saturations, pressure, and temperature. Our choices, which are representative of thermal reservoir simulation, are given in \rtwo{Appendix} \ref{sec:murho}. 

Since the saturations satisfy the constraint \eqref{eq:satconstraint}, one of the saturation variables is explicitly replaced
\begin{equation}\label{eq:Sbeta}
S_\beta = 1 - \sum_{\alpha\neq \beta} S_\alpha. 
\end{equation}
Thus, if the following system is for the flow of $n_\mathrm{phases}$ phases, it consists of \rone{$n_\mathrm{phases} + 1$} equations and unknowns.

We are interested in solving the boundary value problem for $p$, $T$, $S_\alpha$ for all $\alpha\neq \beta$ satisfying equations \eqref{eq:consmass}, \eqref{eq:energy3} together with boundary conditions \eqref{eq:massNeumann}, \eqref{eq:massDirichlet}, \eqref{eq:energy4}, where initial conditions for $p$, $T$, and $S_\alpha$ are prescribed, and we explicitly eliminate $S_\beta$ via \eqref{eq:Sbeta}.

The numerical tests in \rone{this paper} describe the flow of two phases: a heavy oil and water, denoted 
$o$ and $w$, respectively.

\subsubsection{Source/sink terms}\label{sec:source}
Efficient recovery techniques for heavy oils include some way of reducing their viscosity, for example by increasing the temperature inside the reservoir. This can be done by injecting hot fluids (for example steam) or by electromagnetic heating \cite{sahni2000electromagnetic}.

We first consider the source/sink terms representing injection and production wells. These can be modelled by using point sources/sinks
\begin{align*}
 f_\alpha(\mathbf x ) =& \sum_i q^i_{\alpha,\mathrm{inj}}(p,T,S_o) \delta (\mathbf x - \mathbf x_\mathrm{inj}^i) \rho_\alpha  (p,T_\mathrm{inj}) \\
         &-\sum_j q^j_{\alpha,\mathrm{prod}}(p,T,S_o) \delta (\mathbf x - \mathbf x_\mathrm{prod}^j) \rho_\alpha (p,T) , \numberthis \label{eq:f}
\end{align*}
\begin{align*}
 f_T(\mathbf x ) =& \sum_\alpha \left( \sum_i q^i_{\alpha,\mathrm{inj}}(p,T,S_o) \delta (\mathbf x - \mathbf x_\mathrm{inj}^i) \rho_\alpha (p,T_\mathrm{inj}) c_{v,\alpha} T_\mathrm{inj} \right. \\
         &\left. -\sum_j q^j_{\alpha,\mathrm{prod}} (p,T,S_o) \delta (\mathbf x - \mathbf x_\mathrm{prod}^j) \rho_\alpha (p,T) c_{v,\alpha} T \right) , \numberthis \label{eq:fT}
\end{align*}
where $\mathbf x_\mathrm{inj}$ and $\mathbf x_\mathrm{prod}$ represent the locations of injection and production wells, respectively, $\delta (\mathbf{x})$ is the Dirac delta function, $q_{\alpha,\mathrm{inj}}^i$ and $q_{\alpha,\mathrm{prod}}^i$ are the injection and production rates of well $i$, respectively, for phase $\alpha$. 

For very simple cases, we use constant rates. In other cases, we consider a version of the commonly used Peaceman well model \cite{chen2009well,peaceman1978interpretation}. Using such a model usually requires solving additional coupled equations for each well. Since we are using point source/sink wells, we do not have to solve any well equations to determine rates. 

For heating source terms, instead of considering an electromagnetic model, we simply use the function 
\begin{equation}
 f_T =  \sum_i U_\mathrm{heater}(p,T) \delta (\mathbf x - \mathbf x_\mathrm{heater}^i) (T_\mathrm{heater} - T),
\end{equation}
where $\mathbf x_\mathrm{heater}$ represent the locations of heaters, $U_\mathrm{heater}$ is the heat transfer coefficient, and $T_\mathrm{heater}$ is the target heating temperature.

Details and \rev{parameters} about well rates and heating sources terms can be found in Section 2.1.5 of \cite{roy2019block}.

\subsection{DG(0) discretization}
\label{sec:fem}

Commercial reservoir simulators typically use first order Finite Volume methods \cite{leveque2002finite} for spatial discretization. These methods are stable, conservative, and easy to implement. 

To have easy access to the different preconditioners provided by the PETSc library \cite{petsc-web-page}, we use the open-source Finite Element software Firedrake \cite{rathgeber2016firedrake} to handle spatial discretization. The user must express the weak form of their PDE system in the Unified Form Language (UFL) \cite{alnaes2014unified}. Therefore, we want to formulate a weak version of the Finite Volume method used in reservoir simulation. To this end, we use a discontinuous Galerkin (DG) method \cite{riviere2008discontinuous} where the facet integrals are approximated as in reservoir simulation using two-point flux approximation. The details of the discretization are given in Section 2.2 of \cite{roy2019block}, but we first provide the necessary notation.

Let $\mathcal{T} = \{E_i, i\in \mathcal{I}\}$ be a given partition of $\Omega$ into open element domains $E_i$ such that the union of their closure is $\overbar \Omega$, where $\mathcal{I}$ is a set of indices. The $E_i$ are called cells. Since we want to emulate the Finite Volume method, we consider a piecewise constant approximation of our solution variables. The approximate solution will lie in the approximation space $\mathcal{V}_h = \mathbb P_\mathrm{DG}^0$ with the basis of cell-wise indicator functions, $\{\phi_i = \mathbbm{1}_{E_i}\mid i\in \mathcal{I} \}$. The DG(0) approximation of a function $u$ is $u_h = \sum_{i\in\mathcal{I}} u_i \phi_i$, with real coefficients $u_i$.

For a given ordering of the indices in $\mathcal{I}$, we denote by $u^+$ and $u^-$ the limit values of $u$ for two cells sharing an edge. On an edge, the jump of $v$ is defined as $[v] = v^+ - v^-$. Denoting cell centers by $h^+$ and $h^-$, we get that $\|[h]\|$ is an edge quantity representing the Euclidean distance between the adjacent cell centers.

\subsubsection{Semidiscrete problem}
We now discretize the system \eqref{eq:consmass}--\eqref{eq:massDirichlet}, \eqref{eq:energy3}, \eqref{eq:energy4} in space using the semidiscrete DG(0) formulation described in \cite{roy2019block}. We consider the discretization for two-phase flow with an oil phase and a water phase, denoted $o$ and $w$, respectively. Here, we set $S_w = 1 - S_o$. Assuming homogeneous Neumann boundary conditions, the variational problem is: find the approximation $(p,T,S_o)\in \mathbb P_\mathrm{DG}^0 \times \mathbb P_\mathrm{DG}^0 \times \mathbb P_\mathrm{DG}^0$ satisfying the conservation of water mass equation: 
\begin{multline}\label{eq:conswater}
\int_\Omega \phi\pder{(\rho_w (1-S_o))}{t} q \diff x \\+  \int_{\Gamma_\mathrm{int}} [q] \left( \{\!\!\{\mathbf{K}\}\!\!\}\frac{\rho_w^\mathrm{up}k_{rw}^\mathrm{up}}{\mu_w^\mathrm{up}}\left(\frac{[p]}{\|[h]\|} - \{\rho_w\} \mathbf{g} \cdot\mathbf n_e\right)\right)\diff S - \int_\Omega f_w q \diff x = 0, 
\end{multline}
the conservation of energy equation:
\begin{multline}\label{eq:consenergy}
 \int_\Omega \phi\pder{(c_{v,w}\rho_w (1-S_o)T)}{t} r \diff x +  \int_\Omega \phi\pder{(c_{v,o}\rho_o S_o T)}{t} r \diff x + \int_\Omega (1-\phi) \rho_r c_r \pder{T}{t} r \diff x \\
 + \sum_{\alpha= o,w} \int_{\Gamma_\mathrm{int}} [r]  \{\!\!\{\mathbf{K}\}\!\!\}k_{r\alpha}^\mathrm{up} c_{v,\alpha}\frac{\rho_\alpha^\mathrm{up}}{\mu_\alpha^\mathrm{up}}T^\mathrm{up}\left(\frac{[p]}{\|[h]\|} - \{\rho_\alpha\} \mathbf{g} \cdot\mathbf n_e\right)\diff S \\
 + \int_{\Gamma_\mathrm{int}} [r] \{\!\!\{k_T\}\!\!\}\frac{[T]}{\|[h]\|}\diff S- \int_\Omega f_T r \diff x = 0,
\end{multline}
and the conservation of oil mass equation:
\begin{align*}
\int_\Omega \phi\pder{(\rho_o S_o)}{t} s \diff x +  \int_{\Gamma_\mathrm{int}} [s] \left( \{\!\!\{\mathbf{K}\}\!\!\}\frac{\rho_o^\mathrm{up}k_{ro}^\mathrm{up}}{\mu_o^\mathrm{up}}\left(\frac{[p]}{\|[h]\|} - \{\rho_o\} \mathbf{g} \cdot\mathbf n_e\right)\right)\diff S \\- \int_\Omega f_o s \diff x = 0, \numberthis \label{eq:consoil}
\end{align*}
for all $(q,r,s)\in \mathbb P_\mathrm{DG}^0 \times \mathbb P_\mathrm{DG}^0\times \mathbb P_\mathrm{DG}^0$. The unit outward pointing normal of a cell is denoted by $\mathbf n_e$. The single brackets $\{\cdot\}$ denote the average across facets, while the double brackets $\{\!\!\{\cdot\}\!\!\}$ denote the harmonic average across facets. This choice of averaging is standard in reservoir simulation \rtwo{\cite{eymard2000finite, patankar1980numerical}.}

Since the direction of Darcy flow is determined by the pressure gradient and gravity, upwind quantities carried by the velocity $\mathbf{u}_\alpha$\rev{, which, for a facet $e$ shared by $E_1$ and $E_2$ and $\mathbf n_e$ pointing from $E_1$ to $E_2$,} are given by
\begin{equation}
 (u_\alpha)^\mathrm{up} = \begin{cases}
                  u\mid_{E_1}^e    &\text{if  }  \dfrac{[p]}{\|[h]\|} - \{\rho_\alpha\} \mathbf{g} \cdot\mathbf n_e \geq 0,\\
                  u\mid_{E_2}^e    &\text{if  }  \dfrac{[p]}{\|[h]\|} - \{\rho_\alpha\} \mathbf{g} \cdot\mathbf n_e < 0.
                 \end{cases}
\end{equation}
The delta functions in the source/sink terms are approximated on coarse grids as
\begin{equation}
 \delta (x) = \begin{cases}
            1/|E_i| &\text{if } x \in E_i, \\
            0 &\text{otherwise}.
            \end{cases}
\end{equation}
For mesh refinement cases, the volume of the approximated delta functions are kept constant so that they will be defined over several cells. 

\subsubsection{Fully discretized problem}\label{sec:fully}
Time discretization is done using the fully implicit Euler method. We define three mappings, which are linear with respect with their last argument. \rone{The mappings are defined by the implicit Euler scheme applied to the left-hand side of the conservation equations \eqref{eq:conswater}-\eqref{eq:consoil}. We denote the conservation of water mass mapping, the conservation of energy mapping, and conservation of oil mapping by $F_w(p, T, S_o; q)$, $F_e(p,T,S_o;r)$, and $F_o(p,T,S_o;s)$, respectively.}
Let 
\begin{equation}\label{eq:F=0}
 F(p,T,S_o;q,r,s) \coloneqq F_w(p,T,S_o;q) +  F_e(p,T,S_o;r) + F_o(p,T,S_o;s),
\end{equation}
which is linear in $q$, $r$, and $s$, but nonlinear in $p$, $T$, and $S_o$. 

At each time-step, given the previous solution $(p^n,T^n, S_o^n) $, we search for \\ $(p^{n+1},T^{n+1},S_o^{n+1})\in\mathbb P_\mathrm{DG}^0\times\mathbb P_\mathrm{DG}^0\times\mathbb P_\mathrm{DG}^0$ such that
\begin{equation}\label{eq:nonlinear}
 F(p^{n+1},T^{n+1},S_o^{n+1};q,r,s) = 0\quad \text{for all } (q,r,s)\in \mathbb P_\mathrm{DG}^0\times\mathbb P_\mathrm{DG}^0\times\mathbb P_\mathrm{DG}^0.
\end{equation}
\rone{Consider the equivalent problem}: find  $(p^{n+1},T^{n+1},S_o^{n+1})\in\mathbb P_\mathrm{DG}^0\times\mathbb P_\mathrm{DG}^0\times\mathbb P_\mathrm{DG}^0$ such that
\begin{equation}\label{eq:nonlinear2}
 F^*(p^{n+1},T^{n+1},S_o^{n+1};q,r,s) = 0\quad \text{for all } (q,r,s)\in \mathbb P_\mathrm{DG}^0\times\mathbb P_\mathrm{DG}^0\times\mathbb P_\mathrm{DG}^0,
\end{equation}
where
\begin{equation}\label{eq:Fstar}
 F^*(p,T,S_o;q,r,s) \coloneqq F_p (p,T,S_o;q) +  F_e(p,T,S_o;r) + F_o(p,T,S_o;s)
\end{equation}
is $F$ where $F_w$ has been replaced by the following ``pressure equation'' mapping:
\begin{equation}\label{eq:weighted}
 F_p (p,T,S_o;q) \coloneqq c_{v,w} F_w(p,T,S_o;q) + c_{v,o} F_o(p,T,S_o;q).
\end{equation}
Recall that $c_{v,w}$ and $c_{v,o}$ are the specific heat coefficients of the water and oil phases, respectively.


\section{Solution algorithms}
\label{sec:precons}

The system \eqref{eq:nonlinear2} can be written as a nonlinear system of equations for the real coefficients $p_i$, $T_i$ and $S_{o,i}$ of the DG(0) functions $p^{n+1}$, $T^{n+1}$, and $S_o^{n+1}$, respectively. Let $x$ be the vector of these real coefficients and \rone{$\mathcal {F}$} the function such that $\mathcal{F}(x) = 0$ is equivalent to \eqref{eq:nonlinear2}. Solving this equation using Newton's method requires the solution of the following linearized system at each iteration $k+1$:
\begin{equation}
 \pder{\mathcal{F}}{x}\mid_{x=x_k} (x_{k+1} - x_k) = -\mathcal{F}(x_k),
\end{equation}
which can be written as a block system of the form
\begin{equation}\label{eq:AXR}
A \delta x = \begin{bmatrix}
A_{pp} & A_{pT} & A_{ps}\\
A_{Tp} & A_{TT} & A_{Ts}\\
A_{sp} & A_{sT} & A_{ss}
\end{bmatrix}
\begin{bmatrix}
\delta p \\ \delta T \\ \delta s
\end{bmatrix}
= \begin{bmatrix}
b_p \\ b_T \\ b_s
\end{bmatrix} = b,
\end{equation}
where $\delta x = x_{k+1} - x_k$, and $b = -\mathcal{F}(x_k)$, in the desired ordering. \rtwo{Looking back at the strong form of the equations \eqref{eq:consmass} \eqref{eq:energy3}, we see that} the blocks of the pressure-temperature submatrix are the discrete versions of 

\begin{equation}\label{eq:App}
 A_{pp} \sim \phi \frac{1}{\Delta t} \sum_\alpha c_{v,\alpha}S_\alpha(\rho_\alpha)_p + \sum_\alpha c_{v,\alpha} \nabla \cdot \left(\rho_\alpha \mathbf{u}_\alpha\right)_p - \sum_\alpha c_{v,\alpha} (f_\alpha)_p,
\end{equation}
\begin{equation}
 A_{pT} \sim \phi \frac{1}{\Delta t} \sum_\alpha c_{v,\alpha} S_\alpha (\rho_\alpha)_T  + \sum_\alpha c_{v,\alpha} \nabla \cdot \left(\rho_\alpha \mathbf{u}_\alpha \right)_T - \sum_\alpha c_{v,\alpha} (f_\alpha)_T,
\end{equation}
\begin{equation}
 A_{Tp} \sim \phi \frac{1}{\Delta t} \sum_\alpha c_{v,\alpha} S_\alpha (\rho_\alpha)_p T + \sum_\alpha \nabla \cdot \left(c_{v,\alpha} T(\rho_\alpha \mathbf{u}_\alpha)_p\right) - (f_T)_p,
\end{equation}
\begin{multline}
  A_{TT} \sim \frac{\phi}{\Delta t} \sum_\alpha c_{v,\alpha} (\rho_\alpha + (\rho_\alpha)_T T)  + (1-\phi) \frac{\rho_r c_r }{\Delta t} + \sum_\alpha \nabla \cdot \left( c_{v,\alpha} \rho_\alpha \mathbf{u}_\alpha\right) \\ + \sum_\alpha \nabla \cdot \left( c_{v,\alpha}  T (\rho_\alpha \mathbf{u}_\alpha)_T\right) - \nabla \cdot (k_T \nabla) - (f_T)_T,
\end{multline}
where 
\begin{equation}
 (\rho_\alpha \mathbf{u}_\alpha)_p = -\frac{\mathbf K k_{r\alpha}}{\mu_\alpha} \left(\rho_\alpha (\nabla - (\rho_\alpha)_p \mathbf{g})+ (\rho_\alpha)_p (\nabla p - \rho_\alpha \mathbf{g} )\right),
\end{equation}
and
\begin{equation}\label{eq:rhouT}
 (\rho_\alpha \mathbf{u}_\alpha)_T = -\mathbf K k_{r\alpha} \left[ \left( \frac{\rho_\alpha}{\mu_\alpha} \right)_T (\nabla p - \rho_\alpha \mathbf{g} ) - \frac{\rho_\alpha}{\mu_\alpha}  (\rho_\alpha)_T \mathbf{g}\right].
\end{equation}
Here, the coefficients in \eqref{eq:App}-\eqref{eq:rhouT} are evaluated at the previous Newton iterate $(p_k, T_k, S_{o,k})$, and $(.)_p$ and $(.)_T$ denote the partial derivatives with respect to $p$ and $T$, respectively.

We seek to solve the resulting linearized system \eqref{eq:AXR} using iterative methods \cite{saad2003iterative}. We will consider the generalized minimal residual method (GMRES) \cite{saad1986gmres}, a widely-used Krylov subspace method suitable for nonsymmetric linear systems.


For \rtwo{complicated} linear systems such as the ones considered in this paper, an efficient preconditioner is required in order to achieve a rapid convergence with a Krylov method. In this section, we will detail different preconditioning techniques that take into account the specific structure of \eqref{eq:AXR}. 
\rtwo{Important ingredients for the preconditioning techniques are Incomplete LU factorization (ILU) \cite{meijerink1977iterative,saad2003iterative} and Algebraic Multigrid (AMG) \cite{ruge1987algebraic,stuben2001introduction}.}

%

\subsubsection*{Two-stage preconditioners}
The methods described in Sections \ref{sec:cpr} and \ref{sec:cptr} are multiplicative two-stage preconditioners, and so we first provide here a generic definition. Let $M_1$ and $M_2$ be two preconditioners approximating the system matrix $A$, for which we know a way of applying their (generally approximate) inverses $M_1^{-1}$, and $M_2^{-1}$. Applying the resulting two-stage preconditioner for the system $Ax=b$ can be done as follows:
\begin{enumerate}
 \item Apply the first preconditioner: $x_1 = M_1^{-1}b$;
 \item Compute the new residual: $b_1 = b - A x_1$;
 \item Apply the second preconditioner: $x = M_2^{-1} b_1 + x_1$.
\end{enumerate}
The action of the two-stage preconditioner can also be written in the explicit form
\begin{equation}
 M^{-1} = M^{-1}_{2} (I - AM_{1}^{-1}) + M^{-1}_1.
\end{equation}

\subsection{Constrained pressure residual (CPR)}
\label{sec:cpr}
The standard preconditioner for multiphase flow in porous media is the Constrained Pressure Residual method (CPR) \cite{wallis1983incomplete}. Initially developed for isothermal cases, CPR is also used in the thermal case where the temperature variable is treated like the saturation variables. The linear system \eqref{eq:AXR} is thus rearranged in the following form:
\begin{equation}\label{eq:AXRps}
A \delta x = \begin{bmatrix}
A_{pp} & A_{ps}\\
A_{sp} & A_{ss}
\end{bmatrix}
\begin{bmatrix}
\delta p  \\ \delta s
\end{bmatrix}
= \begin{bmatrix}
b_p \\ b_s
\end{bmatrix} = b,
\end{equation}
where the temperature is grouped with the saturations.

The CPR preconditioner is a two-stage preconditioner where the preconditioner $M_1$ solves a restricted system for the pressure
\begin{equation}  
    M_1^{-1} \approx \begin{bmatrix}
                            A_{pp}^{-1} & 0 \\
                            0 & 0
                            \end{bmatrix},
\end{equation}
 where $A_{pp}^{-1}$ is generally approximated using an AMG V-cycle. 
 
 If the pressure solution given in the first stage is accurate, only transport problems for the saturations and temperature remain (assuming that heat diffusion is not significant). Therefore, a simple method such as \rone{ILU(0) or ILU(1)} is adequate for the second stage preconditioner $M_2$. 
 
 Applying CPR to the right-hand side $b = \rvec{ b_p , b_s}$ can be done as follows:
 \begin{enumerate}
  \item Solve the pressure subsystem: $A_{pp} x_p = b_p$;\vspace{0.5em}
  \item Compute the new residual: $\tilde{b} = \begin{bmatrix} b_p \\ b_s \end{bmatrix} - \begin{bmatrix}A_{pp} \\ A_{ps} \end{bmatrix} x_p$;
  \vspace{0.5em}
  \item Precondition and correct: 
  $\begin{bmatrix} \delta_p \\ \delta_s \end{bmatrix} = M_2^{-1}  \tilde b + \begin{bmatrix} x_p\\0 \end{bmatrix}.$
 \end{enumerate}
Thus, the update provided by CPR is
\begin{equation}
  \delta = M_2^{-1}\left(I-(A-M_2)\begin{bmatrix}
                             A_{pp}^{-1} & 0 \\
                             0 & 0 
                            \end{bmatrix}\right) b.
\end{equation}
 Note that an approximation of the pressure equation $A_{pp} \delta p + A_{ps} \delta s= b_s$ \rtwo{in \eqref{eq:AXRps}} is performed in the first stage of CPR where the coupling with secondary unknowns $A_{ps}$ is ignored. To reduce the coupling between the pressure equation and the secondary unknowns, decoupling operators are often used with CPR. Applied to the linear system \eqref{eq:AXRps}, decoupling operators have the form
\begin{equation}G=
 \begin{bmatrix}
       I & -D \\
       0 & I 
      \end{bmatrix},
\end{equation}
where $D$ is a block diagonal matrix \rtwo{(with wide rectangular blocks if there are multiple non-pressure unknowns)}. Decoupling operators are often applied to the linear system as a preprocessing step, i.e. a left preconditioner. Alternatively, the decoupling operator $G$ can be used within the first stage of CPR as follows:
\begin{equation}
M_1^{-1} = R_p^\top(R_p G A R_p^\top)^{-1}R_p G,
\end{equation}
where $R_p$ is the pressure restriction operator such that $R_p A R_p^\top = A_{pp}$.

The most common approximations are Quasi-IMPES (QI) and True-IMPES (TI) ~\cite{lacroix2003iterative,lacroix2000iterative,scheichl2003decoupling}. These approximations are $D_\mathrm{QI} =\mathrm{diag} (A_{ps})\, \mathrm{diag} (A_{ss})^{-1}$, $D_\mathrm{TI} = \\ \mathrm{colsum} (A_{ps}) \mathrm{colsum} (A_{ss})^{-1}$, respectively. Here, $\mathrm{diag}(.)$ and $\mathrm{colsum}(.)$ return block diagonal matrices. \rone{The $\mathrm{colsum}(.)$ operator returns the sum of the entries in each column (summing the columns of block matrices independently).} \rtwo{If the matrices are non-square, the $\mathrm{diag}(.)$ and $\mathrm{colsum}(.)$ operators return block diagonal matrices with non-square blocks.}

\rone{When using decoupling operators, the first stage of CPR consists of solving the approximate Schur complement system involving $S_p = A_{pp} - D A_{sp}$ instead of the original pressure block $A_{pp}$
}. While the elliptic-like properties of $A_{pp}$ are clear from the PDEs, it is less so for $S_p$. Properties amenable to the performance of AMG (e.g. M-matrix properties) are nearly guaranteed in simpler multiphase cases \cite{gries2015system}, but compositional and thermal cases have no such guarantee.

In our case, we also perform a weighted sum of the mass equations in \eqref{eq:weighted}. This is equivalent to choosing $D$ as a block diagonal matrix where the coefficients are $c_{v,\alpha}$ for the mass conservation equation of phase $\alpha$ and 0 for the energy conservation equation.

For our test cases, using decoupling operators can reduce the number of GMRES iterations needed. TI is observed to be superior to QI, which sometimes increases the number of GMRES iterations.

\subsection{Constrained pressure-temperature residual (CPTR)}\label{sec:cptr}
We introduce a CPR-like two-stage preconditioner where a pressure-temperature subsystem is solved approximately in the first stage using an extension of the block preconditioner from \cite{roy2019block}. We call this method Constrained Pressure-Temperature Residual (CPTR).

Let the pressure-temperature submatrix 
\begin{equation}\label{eq:A00}
A_{00} = 
 \begin{bmatrix}
A_{pp} & A_{pT} \\
A_{Tp} & A_{TT} \\
\end{bmatrix},\quad \text{so that} \quad
A = 
 \begin{bmatrix}
A_{00} & A_{0s} \\
A_{s0} & A_{ss} \\
\end{bmatrix}.
\end{equation}
For CPTR, The first stage preconditioner $M_1$ is given by
\begin{equation}
 M_1^{-1} =  \begin{bmatrix}
              A_{00}^{-1} & 0 \\ 0 & 0
             \end{bmatrix},
\end{equation}
where $A_{00}^{-1}$ is an approximation of the action of the inverse of $A_{00}$. 

In our case, we consider two approaches to approximate $A_{00}^{-1}$: an approximation given by a Schur complement factorization detailed in Section \ref{sec:block}, and an unknown-based AMG method discussed in Section \ref{sec:AMGpT}. 

Similarly to CPR, we use ILU for the second stage preconditioner $M_2$.

Let $b_0 = \rvec{b_p, b_T}$. Applying CPTR to the right-hand side $b = \rvec{ b_0, b_s}$ can be done as follows:
 \begin{enumerate}
  \item Solve the pressure-temperature subsystem: $A_{00}  x_0 = b_0$;\vspace{0.5em}
  \item Compute the new residual: $\tilde{b} = \begin{bmatrix} b_0 \\ b_s \end{bmatrix} - \begin{bmatrix}A_{00} \\ A_{0s} \end{bmatrix} x_0$;
  \vspace{0.5em}
  \item Precondition and correct: 
  $\begin{bmatrix} \delta_0 \\ \delta_s \end{bmatrix} = M_2^{-1}  \tilde b + \begin{bmatrix} \rtwo{x_0}\\0 \end{bmatrix}.$
 \end{enumerate}
Thus
\begin{equation}
  \delta = M_2^{-1}\left(I-(A-M_2)\begin{bmatrix}
                             A_{00}^{-1} & 0 \\
                             0 & 0 
                            \end{bmatrix}\right) b.
\end{equation}

Similarly to CPR, an approximation is made in the first stage of CPTR where the saturation coupling block $A_{0s}$ is ignored. A first attempt at reducing this decoupling error is to extend the idea of decoupling operators to the pressure-temperature system. In that case, the first stage of CPTR consists in solving the system with the Schur complement approximation $S_0 = A_{00} - DA_{s0}$, for $D$ a block diagonal matrix of appropriate dimensions. The extensions of Quasi-IMPES and True-IMPES are $D_\mathrm{QI}=\mathrm{diag}(A_{0s}) \mathrm{diag}(A_{ss})^{-1}$, and $D_\mathrm{TI}=\mathrm{colsum}(A_{0s}) \mathrm{colsum}(A_{ss})^{-1}$, respectively.


\subsection{Block preconditioner for the pressure-temperature subsystem}
\label{sec:block}
In this section, we present a choice for the pressure-temperature solver $A_{00}^{-1}$. Consider the following decomposition of the pressure-temperature submatrix:
 \begin{equation}\label{eq:Afact}
 A_{00}= 
  \begin{bmatrix}
   I & 0 \\
   A_{Tp} A_{pp}^{-1} & I
  \end{bmatrix}
  \begin{bmatrix}
   A_{pp} & 0 \\
   0 & S_T
  \end{bmatrix}
  \begin{bmatrix}
   I & A_{pp}^{-1} A_{pT} \\
   0 & I
  \end{bmatrix},
 \end{equation}
 where $S_T = A_{TT} - A_{Tp} A_{pp}^{-1} A_{pT}$ is the Schur complement. The inverse of $A_{00}$ is given by
 
  \begin{equation}\label{eq:Ainv}
 A_{00}^{-1} = 
  \begin{bmatrix}
   I & -A_{pp}^{-1} A_{pT} \\
   0 & I
  \end{bmatrix}
  \begin{bmatrix}
   A_{pp}^{-1} & 0 \\
   0 & S_T^{-1}
  \end{bmatrix}
  \begin{bmatrix}
   I & 0 \\
   - A_{Tp} A_{pp}^{-1} & I
  \end{bmatrix}.
 \end{equation}
 In general, the Schur complement $S_T$ is a dense matrix, even if $A_{00}$ is sparse. A common preconditioning technique for PDE systems is to use the blocks of a factorization like \eqref{eq:Afact}, but with a sparse approximation of the Schur complement \cite{elman2014finite}. With an appropriate Schur complement approximation $\tilde{S}_T$, our block preconditioner for $A_{00}$ is given by
 
 \begin{equation}\label{eq:Aschur}
    \tilde{A}_{00}^{-1} = \begin{bmatrix}
   I & -A_{pp}^{-1} A_{pT}  \\
  0 & I
  \end{bmatrix}
  \begin{bmatrix}
   A_{pp}^{-1} & 0 \\
   0 & \tilde S_T^{-1}
  \end{bmatrix}
  \begin{bmatrix}
   I & 0\\
   - A_{Tp} A_{pp}^{-1} & I.
  \end{bmatrix}.
\end{equation}
Both $A_{pp}^{-1}$ and $\tilde S_T^{-1}$ are approximated using an AMG V-cycle. 

Applying the block preconditioner to the right-hand side $\rvec{ b_p , b_T}$ can be done as follows:
 \begin{enumerate}
  \item Solve the pressure subsystem: $A_{pp} x_p = b_p$;
  \item Compute the new energy equation residual: $\tilde{b}_T = b_T - A_{Tp} x_p$;
  \item Solve the Schur complement subsystem: $\tilde{S}_T \delta_T = \tilde{b}_T$;
  \item Compute the new mass equation residual: $\tilde{b}_p = \rone{b_p} - A_{pT} \delta_T$;
  \item Solve the pressure subsystem: $A_{pp} \delta_p = \tilde{b}_p$.
 \end{enumerate}
Thus, $\rvec{\delta_p , \delta_T } = \tilde{A}_{00}^{-1} \rvec{ b_p , b_T }$.

The CPTR method with this pressure-temperature solver involves strictly more work than the first stage of the CPR method described in Section \ref{sec:cpr}. The additional work is essentially an additional application of the solver $A_{pp}^{-1}$, as well as the construction and application of the solver $\tilde{S}_T^{-1}$. One could also use block triangular versions of the block preconditioner, i.e. ignoring the left or right factor in \eqref{eq:Aschur}.

\subsubsection{Schur complement approximation} \label{sec:schur}
 The two simple sparse Schur complement approximations $\tilde{S}_{A_{TT}} = A_{TT}$ and $\tilde{S}_{\mathrm{diag}} = A_{TT} - A_{Tp} \mathrm{diag}\left( A_{pp}\right)^{-1} A_{pT}$ were shown to perform poorly in our previous work for the single phase case \cite{roy2019block}. Here, we extend the Schur complement approximation presented in that work to the multiphase case. This Schur complement approximation requires the construction of the pressure equation \eqref{eq:weighted}.

 We first consider the steady-state problem for pressure and temperature: find $p$, $T$ such that
\begin{equation}\label{eq:steady1}
\sum_\alpha c_{v,\alpha} \nabla \cdot \left(\rho_\alpha \mathbf u_\alpha \right) = 0 \quad \text{in } \Omega.
 \end{equation}
 \begin{equation}\label{eq:steady2}
 \sum_\alpha \nabla \cdot \left(\rho_\alpha c_{v,\alpha} T \mathbf u_\alpha \right) - \nabla \cdot (\mathbf{k}_T \nabla T) = 0 \quad  \text{in } \Omega,
 \end{equation}
where $\mathbf{u}_\alpha$ is given by \eqref{eq:darcy}, and we have homogeneous Neumann boundary conditions. Note that \eqref{eq:steady1} is the pressure equation obtained from the weighting \eqref{eq:weighted}. Here, we consider the linearized pressure-temperature sub-system in a infinite-dimensional setting. Applying a Newton method to \eqref{eq:steady1}-\eqref{eq:steady2}, we obtain a block matrix of the form \eqref{eq:A00} where the blocks are
 \begin{equation}\label{eq:app}
  A_{pp} = \sum_\alpha c_{v,\alpha} \nabla \cdot \left(\rho_\alpha \mathbf u_\alpha \right)_p,
\qquad
  A_{pT} = \sum_\alpha c_{v,\alpha} \nabla \cdot \left(\rho_\alpha  \mathbf u_\alpha \right)_T,
 \end{equation}
 \begin{equation}\label{eq:atp}
  A_{Tp} =  \sum_\alpha \nabla \cdot \left(\rho_\alpha c_{v,\alpha} T \mathbf u_\alpha \right)_p = \sum_\alpha c_{v,\alpha}  \nabla T \cdot \left(\rho_\alpha \mathbf u_\alpha \right)_p + \sum_\alpha c_{v,\alpha} T \nabla \cdot \left(\rho_\alpha  \mathbf u_\alpha \right)_p ,
 \end{equation}
 \begin{align*}
  A_{TT} =&  \sum_\alpha \nabla \cdot \left(\rho_\alpha c_{v,\alpha} T \mathbf u_\alpha \right)_T - \nabla \cdot (\mathbf{k}_T \nabla ) \\
        =& \sum_\alpha c_{v,\alpha} \nabla \cdot \left(\rho_\alpha  \mathbf u_\alpha \right) + \sum_\alpha c_{v,\alpha} \nabla T \cdot \left(\rho_\alpha  \mathbf u_\alpha \right)_T \\ &+ \sum_\alpha c_{v,\alpha} T \nabla \cdot\left(\rho_\alpha  \mathbf u_\alpha \right)_T- \nabla \cdot (\mathbf{k}_T \nabla ) \numberthis \label{eq:att},
 \end{align*}
where the product rule was used for the divergence operator in \eqref{eq:atp} and \eqref{eq:att}. Then, the second term of the Schur complement (which corresponds to the Poincar\'e-Steklov operator in the infinite-dimensional setting) becomes
 \begin{align*}
 A_{Tp} A_{pp}^{-1} A_{pT}
 =& \left(\sum_\alpha c_{v,\alpha}  \nabla T \cdot \left(\rho_\alpha \mathbf u_\alpha \right)_p + \sum_\alpha c_{v,\alpha} T \nabla \cdot \left(\rho_\alpha  \mathbf u_\alpha \right)_p \right) \\
 &\left(\sum_\alpha c_{v,\alpha} \nabla \cdot \left(\rho_\alpha \mathbf u_\alpha \right)_p\right )^{-1}  \sum_\alpha c_{v,\alpha} \nabla \cdot \left(\rho_\alpha  \mathbf u_\alpha \right)_T \\
 =& T\sum_\alpha c_{v,\alpha} \nabla \cdot \left(\rho_\alpha  \mathbf u_\alpha \right)_T + \sum_\alpha c_{v,\alpha}  \nabla T \cdot \left(\rho_\alpha \mathbf u_\alpha \right)_p \\ 
 &\left(\sum_\alpha c_{v,\alpha} \nabla \cdot \left(\rho_\alpha \mathbf u_\alpha \right)_p\right )^{-1}  \sum_\alpha c_{v,\alpha} \nabla \cdot \left(\rho_\alpha  \mathbf u_\alpha \right)_T \\
  =& T\sum_\alpha c_{v,\alpha} \nabla \cdot \left(\rho_\alpha  \mathbf u_\alpha \right)_T +   \nabla T \cdot \sum_\alpha c_{v,\alpha}\left(\rho_\alpha \mathbf u_\alpha \right)_p \\
  &\left( \nabla \cdot \sum_\alpha c_{v,\alpha}\left(\rho_\alpha \mathbf u_\alpha \right)_p\right )^{-1}   \nabla \cdot\sum_\alpha c_{v,\alpha} \left(\rho_\alpha  \mathbf u_\alpha \right)_T. \numberthis \label{eq:secondterm}\\
\end{align*}
 Taking the difference between \eqref{eq:att} and \eqref{eq:secondterm}, we obtain the Schur complement, $A_{TT} - A_{Tp} A_{pp}^{-1} A_{pT}$. We notice that the term $T\sum_\alpha c_{v,\alpha} \nabla \cdot \left(\rho_\alpha  \mathbf u_\alpha \right)_T$ appears in both \eqref{eq:att} and \eqref{eq:secondterm} and thus cancels. This is due to the construction of the pressure equation \eqref{eq:weighted}. We are left with
\begin{multline}\label{eq:fullschur}
 S_T = \sum_\alpha c_{v,\alpha} \nabla \cdot \left(\rho_\alpha  \mathbf u_\alpha \right) - \nabla \cdot (\mathbf{k}_T \nabla ) +\sum_\alpha c_{v,\alpha} \nabla T \cdot \left(\rho_\alpha  \mathbf u_\alpha \right)_T  \\ - \nabla T \cdot \sum_\alpha c_{v,\alpha}\left(\rho_\alpha \mathbf u_\alpha \right)_p \left( \nabla \cdot \sum_\alpha c_{v,\alpha}\left(\rho_\alpha \mathbf u_\alpha \right)_p\right )^{-1}   \nabla \cdot\sum_\alpha c_{v,\alpha} \left(\rho_\alpha  \mathbf u_\alpha \right)_T.
\end{multline}
Since one of the terms coming from the linearization has canceled, we investigate if it is possible that the last two terms in \eqref{eq:fullschur} can cancel under reasonable assumptions. Consider the following operator taken from the last term of \eqref{eq:fullschur}:
\begin{equation}\label{eq:operator}
 \sum_\alpha c_{v,\alpha}\left(\rho_\alpha \mathbf u_\alpha \right)_p \left( \nabla \cdot \sum_\alpha c_{v,\alpha}\left(\rho_\alpha \mathbf u_\alpha \right)_p\right )^{-1}   \nabla \cdot.
\end{equation}
We now make the same approximations as in \cite{roy2019block}, which \rone{assumes} that $\rho_\alpha$ is constant with respect to $p$, and that we use a two-point flux approximation for the facet integrals. This allows us to replace the operator \eqref{eq:operator} by the identity, giving us the following simplified and more practical Schur complement approximation:
\begin{equation}
 \tilde{S}_T = \sum_\alpha c_{v,\alpha} \nabla \cdot \left(\rho_\alpha  \mathbf u_\alpha \right)_T - \nabla \cdot (\mathbf{k}_T \nabla ).
\end{equation}
While the assumption that $\rho_\alpha$ is close to being a constant with respect to $p$ is appropriate for liquid oil and water, it may be less applicable in the case of gases. However, early investigations reveal that increasing the compressibility of the fluid does not reduce the effectiveness of the preconditioner. 

The extension to the time-dependent case and the addition of source/sink terms does not differ from the single phase case, and so we refer the reader to \cite{roy2019block} for more details. In this case, the source/sink terms related to production wells satisfy $f_{T,\mathrm{prod}}= \sum_\alpha c_{v,\alpha} T f_{\alpha,\mathrm{prod}}$, while those for injection wells satisfy $f_{T,\mathrm{inj}} = \sum_\alpha c_{v,\alpha} T_\mathrm{inj} f_{\alpha,\mathrm{inj}}$. Also, the source term for heaters is given by $f_T = U (T_\mathrm{heater} - T) D_\mathrm{heaters}$, where $D_\mathrm{heaters}$ is the sum of delta functions for the locations of heaters. The Schur complement approximation is given by
\begin{multline}\label{eq:schurapproxfull}
  \tilde{S}_{T} = \phi \sum_\alpha \frac{c_{v,\alpha} \rho_\alpha}{\Delta t} + (1-\phi) \frac{\rho_r c_r }{\Delta t} + \sum_\alpha c_{v,\alpha} \nabla \cdot \left(\rho_\alpha  \mathbf u_\alpha \right)_T - \nabla \cdot (k_T \nabla) \\
                 + UD_\mathrm{heaters} -\sum_\alpha c_{v,\alpha} f_{\alpha,\mathrm{prod}}.
\end{multline}
The discretized version of this operator can be obtained \rone{from the conservation of energy mapping $F_e$ from Section \ref{sec:fully}} by removing the terms that only depend on the previous time-step, and by evaluating the coefficients with the values of the previous Newton iterate. We get the following bilinear operator:
\begin{multline} \label{eq:Se}
 S_e(\delta T, r) \coloneqq  \sum_\alpha \int_\Omega \phi c_{v,\alpha}\frac{S_\alpha\rho_\alpha \delta T}{\Delta t} r \diff x 
 + \int_\Omega (1-\phi)\rho_r c_r \frac{\delta T}{\Delta t} r \diff x\\
+ \sum_\alpha \int_{\Gamma_\mathrm{int}} [r]  \{\!\!\{\mathbf{K}\}\!\!\}(k_{r\alpha})^\mathrm{up} c_{v,\alpha}\frac{(\rho_\alpha)^\mathrm{up}}{(\mu_\alpha)^\mathrm{up}}(\delta T)^\mathrm{up}\left(\frac{[p]}{\|[h]\|} - \{\rho_\alpha\} \mathbf{g} \cdot\mathbf n_e\right)\diff S \\
 + \int_{\Gamma_\mathrm{int}} [r] \{\!\!\{k_T\}\!\!\}\frac{[\delta T]}{\|[h]\|}\diff S+ \int_\Omega \left(-\sum_\alpha c_{v,\alpha} f_{\alpha,\mathrm{prod}} + U D_\mathrm{heaters}\right) \delta T\diff x. 
\end{multline}

\rone{Since this Schur complement approximation relies heavily on the PDE structure, it cannot be combined with decoupling operators from Section \ref{sec:cptr} in an obvious way. The following method can be used as an alternative to the block preconditioner.}

\subsection{AMG for the pressure-temperature subsystem}\label{sec:AMGpT}

AMG methods for systems of PDEs are an obvious alternative for the first stage of CPTR. See \cite{clees2005amg} for a detailed description and a general framework for such methods. Efficient implementations of such AMG methods include BoomerAMG \cite{henson2002boomeramg} and multigrid reduction (MGR) from the hypre library \cite{falgout2002hypre}, and Smoothed Aggregation from the ML package \cite{ml-guide}. For example, BoomerAMG was effective in some diffusion-dominated two-phase flow problems \cite{bui2017algebraic}, and MGR had some success with multiphase flow problems \cite{wang2017multigrid, bui2018algebraic}. Such approaches are more easily generalizable than the Schur complement approach since they usually only rely on algebraic information from the system matrix.

The unknown-based strategy is to treat the variables corresponding to the same unknown separately. For a matrix with an unknown-wise ordering, classical AMG coarsening and interpolation are applied to the individual diagonal blocks. Computationally cheap and easy to implement, the unknown-based approach will perform well if the \rone{cross-coupling between unknowns} is not too strong and the diagonal blocks are amenable for the application of classical AMG.

With a typical Finite Volume scheme, degrees of freedoms associated with the different variables are located at the same points. This means that one could consider using a nodal/point-based AMG method (initially designed for linear elasticity problems \cite{ruge1986amg}). Instead of considering all degrees of freedoms individually, they are grouped in a vector at each ``point''. The construction of the coarse problems uses similar heuristics to the scalar equivalent but applied to the ``point'' vectors and resulting blocks in the system matrix. However, the performance of coarsening strategies for the \rone{nodal} approach is problem-specific. Several strategies are discussed in \cite{clees2005amg}. 

Naively applying a \rone{nodal} or unknown-based AMG method to the full reservoir simulation system does not work in the common case where the equations are hyperbolic with respect to the saturations. Substantial work is needed to construct meaningful coarse problems when applying AMG to the full system \cite{clees2010efficient, gries2014preconditioning}. The proposed methods include unknown-based AMG approaches, and decoupling operators focused on maintaining AMG-friendly properties.

However, the equations are elliptic with respect to pressure, and the energy equation is parabolic with respect to temperature. We will consider an unknown-based AMG for the pressure-temperature subsystem \rone{as an alternative to the Schur complement approach}.

\rone{Since AMG methods for systems of PDEs are purely algebraic, they can be applied to the modified pressure-temperature matrix $S_0$ obtained from decoupling operators mentioned in Section \ref{sec:cptr}. It is unclear, however, if the properties of $S_0$ are guaranteed to be amenable to the application of AMG.}





\subsection{Discussion on multi-stage preconditioners}\label{sec:multistage}
The first stage of CPR is given by the approximate solution of $A_{pp}\widetilde{\delta p} = b_p$, which is an approximation of the pressure equation \rtwo{from \eqref{eq:AXR}}
\begin{equation} 
A_{pp} \delta p  + A_{pT} \delta T + A_{ps} \delta s= b_p. 
\end{equation}
If the coupling blocks $A_{pT}$ and $A_{ps}$ are significant, the approximate pressure $\widetilde{\delta p}$ might not be representative. This could be countered by reducing the coupling between the pressure equation and the non-pressure variables. Similarly, the first stage of CPTR is also currently vulnerable to a strong coupling between both the pressure and energy equations and saturation variables since it ignores both $A_{ps}$ and $A_{Ts}$. We will consider the use of decoupling operators for both CPR and CPTR. However, for CPTR, this cannot be done when using our block preconditioner as the first-stage solver.

A simple way to reduce this coupling without the use of algebraic decoupling operators could be to change the order of the two-stage preconditioner \rone{(ILU on the full system followed by the restricted solve)}. Then, the approximate solution of the fully coupled system (through ILU) acts as a decoupling operator. This idea has been used for example in \cite{bui2017algebraic} for the isothermal situation.

For reordered CPTR, investigations revealed that during GMRES the residual of the energy equation was being reduced faster than those of the pressure equation and oil mass conservation equation.
Therefore, to obtain an adequate solution, the equations must be scaled such that they are of similar magnitude. 

In our case, we multiply the pressure equation by a reference temperature, and the oil mass conservation equation by the product of the reference temperature and a reference specific heat coefficient. Thus, we modify $F^*$ from \eqref{eq:Fstar} to
\begin{equation}\label{eq:scaling}
 F^*(p,T,S_o;q,r,s) \coloneqq T_\mathrm{ref} F_p (p,T,S_o;q) +  F_e(p,T,S_o;r) + c_\mathrm{ref} T_\mathrm{ref}F_o(p,T,S_o;s).
\end{equation}
This scaling has little effect on the number of iterations for CPR and the orignal ordering of CPTR. However, it is necessary for \rone{reordered} CPTR, which then performs similarly to the orignal CPTR. Thus, in the next section, we will only consider versions of CPR and CPTR with the original ordering of the stages. Nonetheless, we will use the scaling \eqref{eq:scaling} by default. 

As mentioned above, for \rone{reordered} CPTR, the residual for the energy equation is being reduced faster than the residuals of the mass equations. This indicates that decoupling the energy equation (in this case using ILU) improves the solution of the restricted pressure-temperature subsystem. However, the residual of the other equations is not reduced as fast for this method. We will investigate numerically in the next section if the use of decoupling operators can reduce the number of iterations for CPTR.

Finally, extending CPR and CPTR to more stages does not appear to have a significant effect on the number of iterations. For example, one could consider a three-stage preconditioner where ILU is used before and after a pressure (or pressure-temperature) restricted solution. Even in the cases where this results in fewer iterations, the additional step does not appear to be worth the additional cost. This is also true of the case where a pressure (-temperature) restricted solve is done before and after ILU.


\section{Numerical results}
\label{sec:results}

We perform numerical experiments to evaluate the algorithmic performance of CPR and CPTR. The code used for these experiments can be found on GitHub\footnote{\url{https://github.com/tlroy/thermalporous}}.
This code is implemented using the open-source Finite Element software Firedrake \cite{kirby2018solver,rathgeber2016firedrake}. For solvers, Firedrake interfaces with the PETSc library \cite{petsc-web-page}, allowing for efficient and parallel computations. Our Schur complement approximation and the two-stage preconditioners are implemented through Firedrake's Python interface, including proof-of-concept implementations for the decoupling operators. Simple but less efficient versions of the two-stage preconditioners can also be implemented by simply providing PETSc options.

The first stage of CPTR is given by our Schur complement approximation \eqref{eq:Se} or by an unknown-based AMG. Unless stated otherwise, for both CPR and CPTR, an AMG V-cycle is used to approximate the inverse for the pressure block $A_{pp}$, the approximate Schur complement, and the temperature block $A_{TT}$. For both scalar and unknown-based AMG, BoomerAMG \cite{henson2002boomeramg} from the hypre library \cite{falgout2002hypre} is used with default parameters, i.e. a symmetric-SOR/Jacobi relaxation scheme (one sweep up, one sweep down), Falgout coarsening,
classical Ruge-St\"uben interpolation, and Gaussian Elimination as the coarse grid solver. BoomerAMG is very efficient in parallel. \rone{When direct LU factorization is used instead of AMG, we use PETSc's default implementation.} For the second stage of CPR and CPTR, ILU(k) is used as provided in PETSc. Unless stated that we use ILU(1), ILU(0) is our default. The parallel version of ILU is block Jacobi where ILU(k) is used for each block (the blocks are determined when Firedrake assigns different parts of the domain to different processors).

For the nonlinear solver, we use Newton's method with line search, and for the linear solver, right-preconditioned GMRES \cite{saad1986gmres}. For Newton's method, the convergence tolerance for both the relative function norm and relative step size norm is set to $10^{-8}$. For GMRES, the convergence tolerance for the relative residual norm is set to $10^{-8}$.

For all cases, we consider oil and water with densities and viscosities as described in \rtwo{Appendix} \ref{sec:murho}. The other physical parameters are shown in Table \ref{tab:parameters}. These parameters are representative of those used in commercial reservoir simulators. The pressures are taken from the SPE10 test case \cite{christie2001tenth}. We also only consider homogeneous Neumann boundary conditions.

\begin{table}[hbt!]
\centering
 \caption{Physical parameters for test cases} \label{tab:parameters}
 \begin{tabular}{l l}
 \hline
 Initial pressure       & 4.1369$\times 10^7$ Pa        \\ 
 Conductivity of oil    & 0.15 W\,m$^{-1}$\,K$^{-1}$        \\ 
 Conductivity of water  & 0.6005638 W\,m$^{-1}$\,K$^{-1}$   \\ 
 Conductivity of rock   & 1.7295772056 W\,m$^{-1}$\,K$^{-1}$\\ 
 Specific heat of oil   & 2093.4 J\,K$^{-1}$\,kg$^{-1}$     \\ 
 Specific heat of water & 4181.3 J\,K$^{-1}$\,kg$^{-1}$     \\ 
 Specific heat of rock  & 920 J\,K$^{-1}$\,kg$^{-1}$        \\ 
 Density of rock        & 2650 kg\,m$^{-1}$             \\ 
 \hline
 \end{tabular}
\end{table}

For the scaling parameters in \eqref{eq:scaling}, we choose $T_\mathrm{ref} = T_\mathrm{prod}$ and $c_\mathrm{ref} = S_o c_{v,o} + (1-S_o) c_{v,w}$, where $S_o$ is the initial oil saturation in the reservoir, which is uniform in all the test cases here. For the preconditioners we consider, the number of iterations is not very sensitive to the choice of scaling parameters. 

\subsubsection*{Computational cost} We note that for the first stage of CPTR, applying the block preconditioner from Section \ref{sec:block} should be slightly more than three times the cost of applying the first stage of CPR, while the unknown-based AMG should be around two times the cost. The relative cost of the second stage will depend on things such as the size and dimension of the problem, and the number of processors used. \rme{Looking} at the CPU timings for the serial case given in Table \ref{tab:paralW} (a problem with 52,728 degrees of freedom), applying the first stage of CPR takes an average of 0.0088s, while the first stage of CPTR with the block preconditioner takes an average of 0.033s, and 0.017s for unknown-based AMG. Applying the second stage of CPR/CPTR with ILU(0) takes an average of 0.02s, and 0.034s with ILU(1). These timings were performed on a Lenovo ThinkCentre M920q with Intel(R) Core(TM) i5-8500T CPU @ 2.10GHz. \rtwo{In terms of memory requirements, the first stage of CPTR requires two different AMG hierarchies. How this compares to the memory requirements of the second stage will depend on the size of the blocks and fill-in used for ILU.} Note that \rone{since we implement our Schur complement approximation in Python,} the implementation of the block preconditioner is not as optimized as \rone{that of unknown-based AMG in BoomerAMG}. Additionally, the full factorization \eqref{eq:schurapproxfull} used here may not be needed in all cases. Dropping one of the factors makes the cost of applying the block preconditioner slightly more than two times the cost of the first stage of CPR, or slightly more than unknown-based AMG.

\subsection{Homogeneous test cases}\label{sec:homogeneous}
We compare the performance of several methods for problems with homogeneous permeability. We begin by a brief mesh refinement study for 2D cases. Then, we look at both weak and strong scaling for 3D cases. For weak scaling, the problem size is increased proportionally with the number of processors, while for strong scaling, the problem size remains fixed. For weak scaling, the problem size is increased via uniform mesh refinement.

For traditional reservoir simulation, strong scaling is often more relevant than weak scaling. Indeed, reservoir models are typically given with geological properties on a (usually rather coarse) fixed grid, for example, the SPE10 test case. Another good reason to study strong scaling is to isolate the effects of parallelization from the effects of mesh refinement observed for weak scaling. Nonetheless, the effects of mesh refinement are very relevant in the context of reservoir simulation, for example around wells or other features of interest.

\subsubsection{Mesh refinement study}  \label{sec:refine}
We perform a mesh refinement study for 2D test cases. The domain is $50\times 50$ meters with an $N \times N$ grid. Starting at $N=20$, we double $N$ until we reach $N=\rme{640}$. Since this is a 2D case, we do not include gravity. For the isotropic cases, the permeability is $3\times10^{-13}$~m$^2$ and the porosity is 0.2. For the anisotropic case, the permeability in the $x$-direction is $3\times10^{-13}$~m$^2$ and in the $y$-direction, $3\times10^{-10}$~m$^2$. 

We consider three cases: one with heaters, and two with wells (isotropic and anisotropic). For each case, we have 6 source/sink terms with 3 located on \rone{two opposing sides} of the square domain. In the heating case, these terms are all heaters, while for the well cases, one side has injection wells and the other, production wells. For the well cases, the production/injection rates are constant at $q=3\times10^{-7}$~m$^3$\,s$^{-1}$, and water is injected at a temperature of 373.15~K. For the heater case, we have a heating temperature of 373.15~K. The initial oil saturation for both cases is set to $S_o=0.9$. We then take two time-steps of 10 days. 

For CPTR, we will compare the use of the following first stage solvers: block preconditioner (\rtwo{CPTR-Block(AMG)}), unknown-wise AMG (\rtwo{CPTR-uAMG}), unknown-wise AMG with the True-IMPES decoupling operator (\rtwo{CPTR-uAMG(TI)}). ILU(0) is used in the second stage. For CPR, we will compare the use of non-decoupled CPR with ILU(0) (\rtwo{CPR-AMG}) and with ILU(1) (\rtwo{CPR-AMG-ILU(1)}), as well as CPR with True-IMPES and ILU(0) (\rtwo{CPR-AMG(TI)}).


\begin{table}[!htb] 
\centering
\caption{\rone{Mesh Refinement:} 2D Isotropic Well case. Average linear iterations per nonlinear iteration \rboth{(computational time (s))}.}\label{tab:meshW}
 \begin{tabular}{|c c c c c c c|}
 \hline
 method/$N$    & 20 & 40 & 80 & 160 & 320 & 640 \\
 \hline
\rtwo{CPR-AMG}            & 5.29 & 6.29  & 9.14 \rboth{(2.3)} & 15.5 \rboth{(9)} & 31.5 \rboth{(67)} & \rme{65.7} \rboth{(580)}  \\
\rtwo{CPR-AMG(TI)}        & 4.86 & 5.43 & 8.29 & 14.7 & 30.7                       & \rme{65.1}  \\
\rtwo{CPR-AMG-ILU(1)}     & 5.29 & 5.43  & 6.43 \rboth{(2.3)} & 10 \rboth{(9)}  & 20 \rboth{(62)}      & \rme{40.6} \rboth{(472)}  \\
\rtwo{CPTR-Block(AMG)}    & 5.29 & 5.29  & 5.43 \rboth{(2.9)} & 5.64 \rboth{(9.5)} & 6.53 \rboth{(55)} & \rme{7.7} \rboth{(310)}  \\
\rtwo{CPTR-uAMG}          & 5.29 & 5.29  & 6.29 \rboth{(1.4)} & 6.82 \rboth{(6.8)} & 8.24 \rboth{(44)} & \rme{10} \rboth{(258)}  \\
\rtwo{CPTR-uAMG(TI)}      & 4.29 & 5.14 & 5.29 & 5.82 & 7.35                       & \rme{9.26}  \\
 \hline
 \end{tabular}
\end{table}


\begin{table}[!htb] 
\centering
\caption{\rone{Mesh Refinement:} 2D Heater case. Average linear iterations per nonlinear iteration \rboth{(computational time (s))}.}\label{tab:meshH}
 \begin{tabular}{|c c c c c c c|}
 \hline
 method/$N$ & 20 & 40 & 80 & 160 & 320 & 640 \\
 \hline
\rtwo{CPR-AMG}            & \rme{5.33} & \rme{5.5}  & \rme{8.5} \rboth{(1.6)} & \rme{14} \rboth{(4.1)}  & \rme{25.8} \rboth{(18)} & \rme{49.3} \rboth{(81)}\\
\rtwo{CPR-AMG(TI)}        & 5.14  & 5.29  & 7.83       & 14.4       & 25          & \rme{48.8}\\
\rtwo{CPR-AMG-ILU(1)}     & \rme{4.67}  & \rme{5.17}  & \rme{5.83} \rboth{(1.6)}& \rme{9} \rboth{(4.2)}& \rme{16} \rboth{(17)} & \rme{30.8} \rboth{(68)}\\
\rtwo{CPTR-Block(AMG)}    & \rme{4.83}  & \rme{5.17}  & \rme{5.17} \rboth{(2)} & \rme{4.6} \rboth{(4.5)}   & \rme{5} \rboth{(15)} & \rme{4.75} \rboth{(45)} \\
\rtwo{CPTR-uAMG}          & \rme{4.67}  & \rme{5}  & \rme{5} \rboth{(1.1)} & \rme{5} \rboth{(3)} & \rme{5} \rboth{(12)} & \rme{4.75} \rboth{(34)}\\
\rtwo{CPTR-uAMG(TI)}      & 4.86 & 4.71 & 4.5  & 4.6  & 4.6 & \rme{4.25} \\
 \hline
 \end{tabular}
\end{table}

\begin{table}[!htb]
\centering
\caption{\rone{Mesh Refinement:} 2D Anisotropic Well case. Average linear iterations per nonlinear iteration \rboth{(computational time (s))}.}\label{tab:anisomeshW}
 \begin{tabular}{|c c c c c c c|}
 \hline
 method/$N$    & 20 & 40 & 80 & 160 & 320 & 640 \\
 \hline
\rtwo{CPR-AMG}            & 5 & 6 	& 9 \rboth{(2.6)}	& 17 \rboth{(10)}	& 32.3 \rboth{(73)} & \rme{67.2} \rboth{(650)} \\
\rtwo{CPR-AMG(TI)}        & 4.29 	& 5.71 	& 8.25 	& 16.1	& 31.8 & \rme{65.3} \\
\rtwo{CPR-AMG-ILU(1)}     & 4.71 	& 5.43	& 6.5 \rboth{(2.6)}	& 11.1 \rboth{(10)}	& 20.4 \rboth{(67)} & \rme{40.9} \rboth{(514)} \\
\rtwo{CPTR-Block(AMG)}    & 4.43 	& 5.14 	& 5.75 \rboth{(3.2)} & 6.92 \rboth{(11)}	& 7.17 \rboth{(60)} & \rme{7.46} \rboth{(319)} \\
\rtwo{CPTR-uAMG}      & 4.86 	& 5.43	 & 6.25 \rboth{(1.7)} & 8.08 \rboth{(8.6)}	& 10.1 \rboth{(55)} & \rme{12.1} \rboth{(333)}\\
\rtwo{CPTR-uAMG(TI)}  & 4.14	& 5	& 5.5	& 7.33	& 9.44 & \rme{11.6} \\
 \hline
 \end{tabular}
\end{table}

For each case, we calculate the average number of linear iterations per nonlinear iteration \rboth{and computational time, obtained with Intel(R) Xeon(R) Gold 6140 CPU @ 2.30GHz. We do not report timings for the methods using TI decoupling since the current implementation is not efficient (TI is very cheap if implemented appropriately). The results} are shown for the isotropic well case, the heater case, and the anisotropic well case in Tables \ref{tab:meshW}, \ref{tab:meshH}, and \ref{tab:anisomeshW}, respectively. In all cases, we observe that for smaller values of $N$, CPR and CPTR perform similarly. However, as we refine the mesh, the number of iterations for CPR increases significantly, while it does not increase much for CPTR, or even decreases in the heater case. In the well cases, we observe that this increase is least significant for \rtwo{CPTR-Block(AMG)}. In the heater case, \rtwo{CPTR-uAMG(TI)} has the smallest number of iterations. \rme{Note that the number of Newton steps is lower for the heater case and does not increase much as the mesh is refined.}

\rboth{It is important to note that a different variable ordering is used for CPTR-uAMG to allow the use of uAMG on the pressure-temperature system. This results in smaller computational times compared to the ordering used for CPR and CPTR-Block. This is not due to uAMG since we observe the same trend when LU is used as the first stage solver for CPTR. Nevertheless, we observe that CPTR-Block(AMG) is only faster than CPR-AMG and CPR-AMG-ILU(1) for $N=320, 640$. For both well cases at $N=640$, the cost per iteration of CPTR-Block(AMG), CPTR-uAMG, and CPR-AMG-ILU(1) is around 4.5, 3 and \rev{1.3} times more than CPR-AMG, respectively.}

CPR treats temperature like a saturation, essentially assuming that heat is simply being transported by the fluid flow. When this is the case, ILU is sufficient. Since heat diffusion becomes more noticeable as the mesh is refined, the CPR strategy does not hold for finer meshes, and ILU is not enough. In contrast, the treatment of temperature for CPTR does not depend on ILU, but rather AMG. Hence, CPTR tackles heat diffusion appropriately on fine meshes.

\subsubsection{Weak scaling}
For weak scaling, the performance of the methods is evaluated as we increase the number of processors and problem size simultaneously. The domain is a $50\times 50\times 50$ meters box, discretized with an $N\times N \times N$ grid. The permeability is constant at $3\times10^{-13}$m$^2$ and the porosity is 0.2. For the number of processors 1, 2, 4, 8, and 16, we have $N$ = 26, 33, 41, 52, 65, respectively. This results in around $50,000$ degrees of freedom per processor. 

For the cases here, \rone{CPR performs similarly with or without True-IMPES, and the different variants of CPTR also perform similarly}. We therefore only show results for CPR-AMG and \rone{\rtwo{CPTR-Block(AMG)}}. \rev{Computational timings for similar test cases are shown in \cite{roy2019preconditioning}, where the different versions of CPTR outperform CPR for large enough problems.}

We consider \rboth{two cases with wells}. For each case, we have 21 source/sink terms both near the top and bottom of the domain. The top terms are injection wells and the bottom terms, production wells. The production/injection rates are constant at $q=10^{-7}$~m$^3$\,s$^{-1}$. Water is injected at a temperature of 373.15~K.

The initial saturation \rboth{is set to $S_o=0.99$ or $S_o = 1$}. For each case, we take 5 time-steps and look at the average number of linear iterations per nonlinear iteration. The size of the time-steps is $\Delta t=10$ days \rboth{for the case with $S_o =0.99$, and $\Delta t=4$ days for the case with $S_o=1$.}


\begin{table}[!htb] \centering
\caption{Weak scaling: 3D Well case, $S_o = 0.99$. Average linear iterations per nonlinear iteration.}\label{tab:paralW099}
 \begin{tabular}{|c c c c c c|}
 \hline
 Method/Num. proc. & 1 & 2 & 4 & 8 & 16 \\
 \hline
\rtwo{CPR-AMG}        & 6    & 7.78 & 8.31 & 9.54 & 12.4 \\
\rtwo{CPTR-Block(AMG)}       & 5.58 & 5.67 & 5.85 & 6.08 & 6.55 \\
\rtwo{CPR-AMG-ILU(1)} & 5.58 & 7.33 & 8    & 9.08 & 11.8 \\
\hline
\end{tabular}
\end{table}

\rboth{The results for the case with $S_o=0.99$ are shown in Table \ref{tab:paralW099}.
We observe that as we increase the number of processors from 1 to 16, the number of iterations for both versions of CPR nearly doubles, while the number of iterations for CPTR increases by fewer than 20\%.}
Note that for both heating and well cases, different starting constant saturations provide similar scaling results, except when $S_o$ is greater than $0.999$ for the well case.

\begin{table}[!htb] \centering
\caption{Weak scaling: 3D Well case, $S_o = 1$. Average linear iterations per nonlinear iteration.}\label{tab:paralW}
 \begin{tabular}{|c c c c c c|}
 \hline
 Method/Num. proc. & 1 & 2 & 4 & 8 & 16 \\
 \hline
\rtwo{CPR-AMG}        & 14.1 & 16   & 17.3  & 20.6 & 23.5 \\
\rtwo{CPTR-Block(AMG)}       & 14.1 & 15.9 & 17.3 & 20.5 & 23.4 \\
\rtwo{CPR-AMG-ILU(1)} & 11.7 & 14   & 15.3  & 17.7 & 19.3 \\
\hline
\end{tabular}
\end{table}

As shown in Table \ref{tab:paralW}, the results are significantly different for the well case with $S_o=1$. Indeed, the number of iterations for all methods increases by the same factor of $2/3$. We observe that CPTR performs essentially the same as CPR (but at an additional cost). This could indicate that the temperature solution given by the first stage of CPTR is either wrong or trivial. Additionally, CPTR scales better in \rboth{Table \ref{tab:paralW099}}, which suggests again that there is a large decoupling error in the first stage of CPTR. As discussed in Section \ref{sec:multistage}, we suspect this is due to a strong coupling between the temperature and saturation. Indeed, the energy equation solved in the first stage of CPTR uses $S_o$ from the previous iteration, which may be a poor approximation when injecting hot water in a reservoir filled with cold oil.

Furthermore, the well case with $S_o=1$ is very similar to some single phase test cases in \cite{roy2019block}, except that oil is injected instead of water. In that paper, the block preconditioner scales significantly better than CPR in both weak and strong scaling. This is a clear indication that the issue with CPTR is related to the coupling of the saturation variable. \rev{However, using decoupling strategies, either CPTR-uAMG(TI) or reordered CPTR-Block mentioned in Section \ref{sec:multistage}, does not change the poor scaling behaviour illustrated in Table \ref{tab:paralW}.}

\subsubsection{Strong scaling}
We use the same problem as the previous section on the finest mesh ($N=65$). We keep the problem size fixed while increasing the number of processors.


\begin{table}[!htb] \centering
\caption{Strong scaling: 3D Well case, $S_o = 0.99$. Average linear iterations per nonlinear iteration.}\label{tab:strongW099}
 \begin{tabular}{|c c c c c c|}
 \hline
 Method/Num. proc. & 1 & 2 & 4 & 8 & 16 \\
 \hline
\rtwo{CPR-AMG}        & 9.45 & 11.4 & 12   & 12.2 & 12.4 \\
\rtwo{CPTR-Block(AMG)}       & 6.15 & 6.3  & 6.35 & 6.3  & 6.55 \\
\rtwo{CPR-AMG-ILU(1)} & 6.8  & 10.4 & 11   & 11.2 & 11.8 \\
\hline
\end{tabular}
\end{table}

\rboth{In Table \ref{tab:strongW099}, we observe for the case with $S_o=0.99$ that as we increase the number of processors from 1 to 16, the number of iterations increases by around 31\% for \rtwo{CPR-AMG}, 73\% for \rtwo{CPR-AMG-ILU(1)} and 7\% for CPTR.}
\rboth{The bulk of the increase for \rtwo{CPR-AMG-ILU(1)} occurs when going from 1 to 2 processors, around 40\%. This showcases how \rtwo{CPR-AMG-ILU(1)} is much more effective in serial.}

\begin{table}[!htb] \centering
\caption{Strong scaling: 3D Well case, $S_o = 1$. Average linear iterations per nonlinear iteration.}\label{tab:strongW}
 \begin{tabular}{|c c c c c c|}
 \hline
 Method/Num. proc. & 1 & 2 & 4 & 8 & 16 \\
 \hline
\rtwo{CPR-AMG}        & 22   & 22.8 & 22.4 & 22.6 & 23.5 \\
\rtwo{CPTR-Block(AMG)}       & 22.1 & 22.8 & 22.3 & 22.5 & 23.4 \\
\rtwo{CPR-AMG-ILU(1)} & 16.6 & 18   & 18.1 & 18.8 & 19.3 \\
\hline
\end{tabular}
\end{table}

In Table \ref{tab:strongW}, we show the results for the well case for $S_o=1$. We observe a similar increase in the number of iterations for \rtwo{CPR-AMG} and CPTR of 7\% and 6\%, respectively. For \rtwo{CPR-AMG-ILU(1)}, the increase in iterations is 17\%, and half of this increase happens when going from 1 to 2 processors.

Both CPR and CPTR use Block ILU, which becomes weaker as the number of processors increases. In cases other than $S_o=1$, CPTR scales much better than CPR. For those cases, the first stage of CPTR gives an accurate temperature update, while CPR relies on \rone{Block} ILU for its temperature update. 

\subsubsection{Synthesis}

\rboth{
In CPR, temperature-related effects are only treated in the second stage by ILU. On the other hand, the CPTR strategy tackles these effects in first stage using some type of AMG method, which captures diffusive effects. As the mesh is refined, heat diffusion becomes more significant, resulting in the poor scaling of CPR, while CPTR exhibits near mesh independence in most well cases, and mesh independence in heater cases. 
}

\rboth{In terms of computational time, CPTR becomes competitive when diffusion is dominant, which in our experiments is when the mesh is sufficiently refined. 
}

\rboth{Block ILU becomes a weaker, albeit cheaper, method as the number of processors increase. Since CPTR does not depend on ILU for the treatment of temperature, it is more robust to a weaker second stage. Note that the authors have observed industrial cases where a strong second stage is needed for CPR's convergence.
}

\rboth{In brief, CPTR is preferable when the second stage \rev{solver} is too weak for CPR, notably for diffusion-dominated flows or when many processors are available.}


\subsection{SPE10 test case}
\label{sec:spe10}

\begin{figure}[hbt]
\centering
 \begin{subfigure}[b]{0.495\textwidth}
 \includegraphics[width=\linewidth]{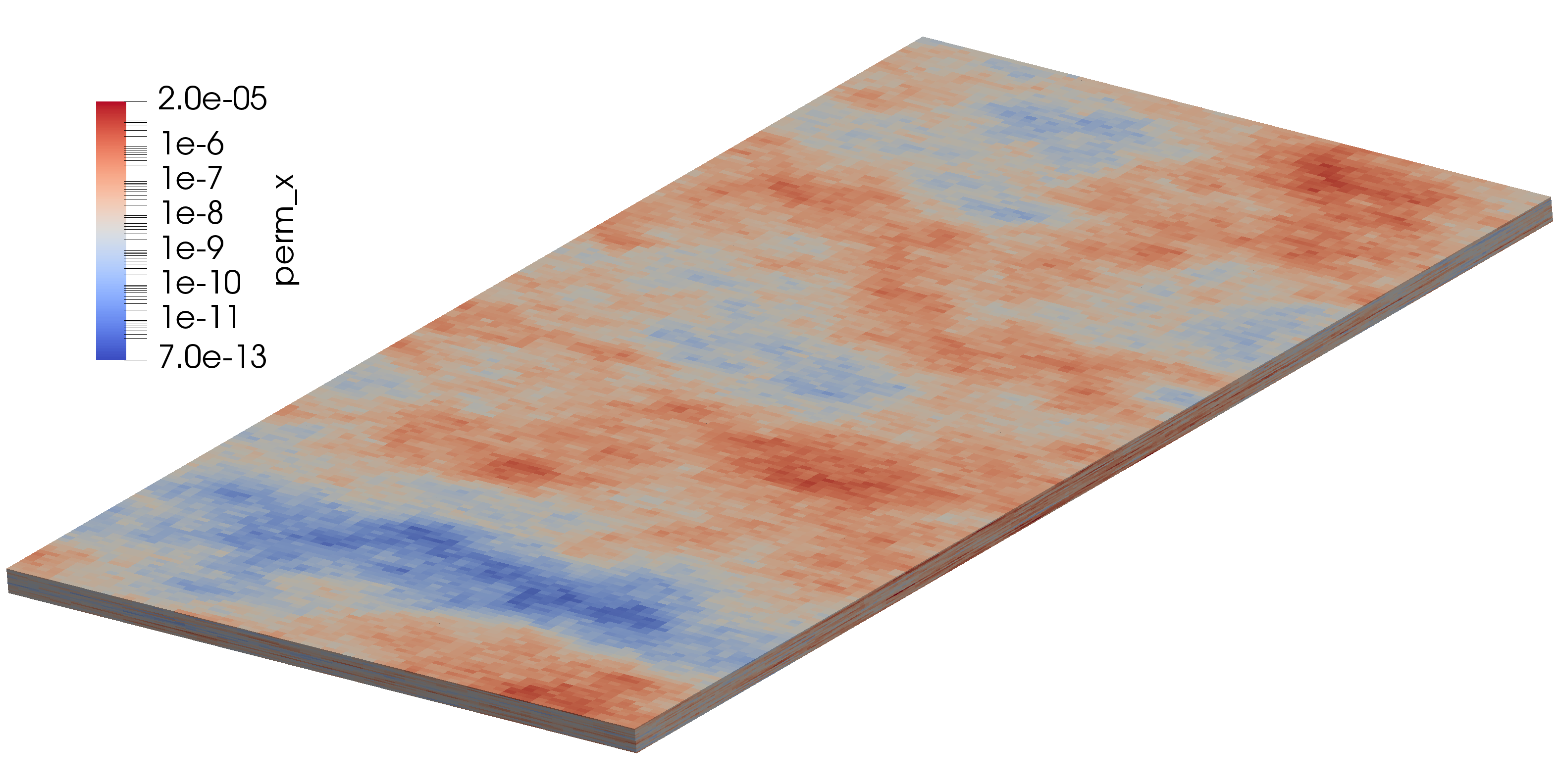}
 \caption{Permeability in $x$ and $y$ directions.}
 \end{subfigure}
 \hfill
 \begin{subfigure}[b]{0.495\textwidth}
 \includegraphics[width=\linewidth]{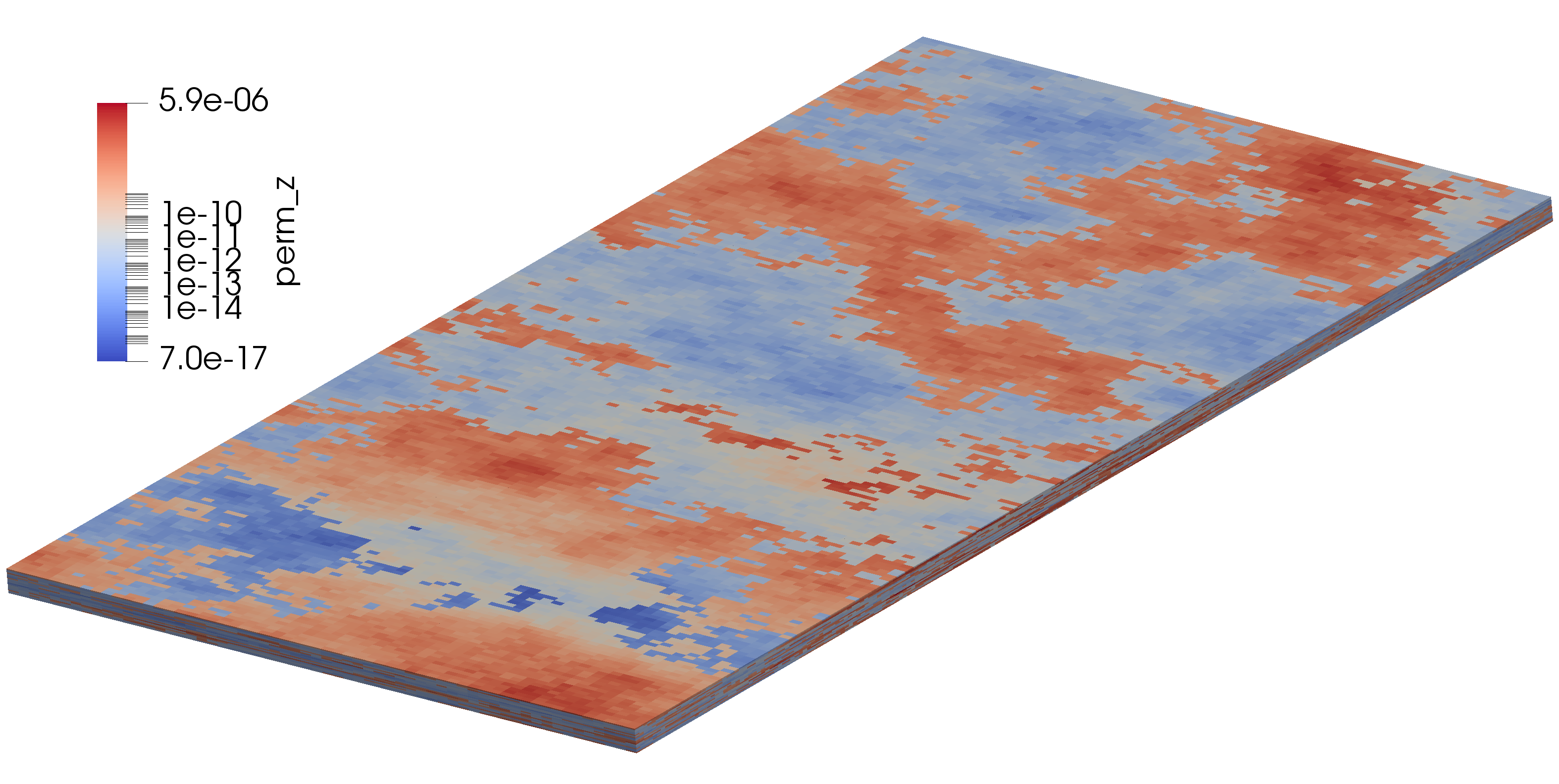}
 \caption{Permeability in $z$ direction}
 \end{subfigure}
 \caption{Log of permeability of the SPE10 test case (mm$^2$).}\label{fig:perms}
\end{figure}

We consider the benchmark problem SPE10 \cite{christie2001tenth}, a highly heterogeneous case, and use its permeability and porosity fields. For the following tests, we consider the top 20 layers of SPE10, such that the domain has dimensions 365.76$\times$670.56$\times$12.192 meters, and the mesh is 60$\times$220$\times$20. The permeability is isotropic in the $x$-$y$ plane, and anisotropic otherwise. The permeability fields are illustrated in Figure \ref{fig:perms}. 

\rboth{We consider a case with wells and heaters (W+H). We position an injection well in the middle of the reservoir, and production wells in each corner. The wells are completed throughout the 20 layers. For the injection and production rates, we use the version of the Peaceman well model described in \cite{roy2019block}. The production wells produce with a bottom-hole pressure of $2.7579\times 10^7$ Pa. The injection well injects hot water with a maximum rate of $q=1.8\times 10^{-3}$~m$^{3}$\,s$^{-1}$. The initial temperature is uniform at 288.706~K, and the injection temperature is 373.15~K. Heater placement is the same as for wells and so is the heating temperature.}

We simulate injection and production for 100 days where the time-steps are chosen with simple heuristics such that Newton's method usually converges in fewer than 10 iterations. For each case, we look at the performance of different methods as we increase the number of processors from 1 to 16. We consider CPR and CPTR where $A_{pp}^{-1}$ and $\tilde{S}_T^{-1}$ are approximated by either AMG V-cyle or direct LU factorization, denoted \rtwo{CPR-AMG and CPTR-Block(AMG), or CPR-LU and CPTR-Block(LU)}, respectively. We include the latter since we believe that there are cases where the AMG V-cycles in the first stage of CPR/CPTR are not as effective as usual. Since CPTR relies on three AMG V-cycles, it is likely more susceptible to any such weaknesses. \rtwo{LU is thus only used as a benchmark for testing the expected lower bound for the iteration count.}

\begin{table}[!htb] \centering
\caption{Strong scaling SPE10 3D wells and heaters case, $S_o = 0.9$. Average linear iterations per nonlinear iteration.}\label{tab:spe10WH}
 \begin{tabular}{|c c c c c c|}
 \hline
 Method/Num. proc. & 1 & 2 & 4 & 8 & 16 \\
 \hline
\rtwo{CPR-AMG}    & 11.5   &11.6   &11.6   &11.6   &11.8 \\
\rtwo{CPTR-Block(AMG)}   & 11.8   &12.2   &12.2   &12.2   &12.3 \\
\rtwo{CPR-LU}     & 10.2   &10.3   &10.2   &10.2   &10.4 \\
\rtwo{CPTR-Block(LU)}    & 10     &9.91   &9.91   &10     &10.1 \\
\hline
\end{tabular}
\end{table}

We take a total of 104 time-steps. The average linear iterations per nonlinear iteration are shown in Table \ref{tab:spe10WH}. We observe that all methods have a small increase in iterations. At the beginning of the simulation, heating is the most significant effect. For the first time-steps, \rtwo{CPTR-Block(AMG)} takes fewer iterations than \rtwo{CPR-AMG}, but the total number of iterations becomes the same at the 30th time-step (around 12 days, \rtwo{results not shown here}). In contrast, \rtwo{CPTR-Block(LU)} has fewer iterations than \rtwo{CPR-LU}. This indicates that single AMG V-cycles are not enough in this case to provide good pressure and temperature solutions in the first stage of \rtwo{CPTR-Block(AMG)}.

 The SPE10 test case was not originally designed for thermal cases, but rather for testing upscaling techniques. Its difficulty lies in its highly heterogeneous permeability field. \rboth{On the other hand, the cases considered in Section \ref{sec:homogeneous} have finer grids and thus diffusion-dominated flows.} \rme{In \cite{cremon2019multi}, versions of CPTR are shown to perform well for SPE10 permeability fields with laboratory length scales.} 

\subsection{Numerical justification of the Schur complement approximation}

We now compare the action of the inverses of the different Schur complement approximations. This is done for cases similar to the single phase cases described in \cite{roy2019block}. The only difference is that we start with $S_o=0.9$, and inject water with a temperature of 373.15~K. In brief, the domain is a 2D slice of SPE10 with a $60\times 120$ grid. For the well case (W), there is one production well and one injection well. These are located in regions of high permeability. For the heater case (H), heater placement is the same as for the well case, and for the well and heater case (W+H), we combine both wells and heaters. For the high permeability cases (h.p.), permeability is increased by a factor of 1,000 to emulate advection-dominated heat flow. 

\begin{table}[!htb]
        \centering
\caption{Condition numbers of the different matrices and Schur complement approximations for various cases ($S_o=0.9$)}
\label{tab:conds}
\renewcommand*{\arraystretch}{1.15}
 \begin{tabular}{|c c c c c c|}
 \hline
 Matrix/Case & H & W & W+H & h.p. W & h.p. W+H \\
 \hline
 $\tilde S_\mathrm{diag}^{-1} S$   & \rtwo{1.043   }& \rtwo{43.11} & \rtwo{122.3} & \rtwo{8.13E+07}& \rtwo{3.64E+07} \\
 $\tilde S_\mathrm{A_{TT}}^{-1} S$ & \rtwo{1.034   }& \rtwo{7.959} & \rtwo{3.050} & \rtwo{6.849   }& \rtwo{382.1} \\
 ${\tilde S_T}^{-1} S$             & \rtwo{1.020   }& \rtwo{1.019} & \rtwo{1.010} & \rtwo{1.556   }& \rtwo{1.95} \\
 $A_{TT}$                          & \rtwo{5.24E+05}& \rtwo{8.802} & \rtwo{1188 } & \rtwo{1093    }& \rtwo{1.60E+05} \\
 $S$                               & \rtwo{5.36E+05}& \rtwo{1.817} & \rtwo{1172 } & \rtwo{241.7   }& \rtwo{1.59E+05} \\
 \hline
 \end{tabular}
\end{table}

For different Schur complement approximations, we look at the condition number of their inverse applied to the full Schur complement \rone{(using a direct solver to invert the matrices)}. While the condition number is not fully indicative of the performance of preconditioned GMRES, it informs us about the quality of the approximations. In Table \ref{tab:conds}, we observe that our custom approximation $\tilde{S}_T$ approximates the Schur complement well, even for the high permeability cases where the simpler approximations struggle.

In terms of the actual number of GMRES iterations, $\tilde{S}_T$ always performs the best for CPTR (results for the single phase case shown in \cite{roy2019block}). This difference is even more significant for harder cases (for example high permeability). In some cases, GMRES does not converge before the prescribed maximum number of iterations when using other approximations. See \cite{roy2019block} for single phase flow examples where the other Schur complements exhibit worse performance for heterogeneous or anisotropic permeability fields.

\subsection{High pressure-temperature cross-coupling} \label{sec:crosscoup}
In Section \ref{sec:refine}, we considered the use of unknown-based AMG (uAMG) and our block preconditioner for the first stage of CPTR. While uAMG usually results in more iterations, it exhibits similar scalability. However, uAMG is not expected to perform well under a strong cross-coupling. Here, we artificially increase the cross-coupling of the pressure-temperature system to compare the robustness of different two-stage preconditioners. Along with uAMG and the block preconditioner, we also show results for CPTR where the first stage is a block diagonal preconditioner using either LU or AMG on the $A_{pp}$ and $A_{TT}$ blocks. To avoid confusion, we call these methods block-diagonal LU \rtwo{(BD(LU))}, and block-diagonal AMG \rtwo{(BD(AMG))}\footnote{\rone{Note that uAMG is different from \rtwo{BD(AMG)} since the coupling blocks remain in the coarse grid correction of uAMG.}}. We also consider CPR with AMG in the first stage.

Recall that we have only considered the flow of two fluids: water and a heavy oil. The viscosity of the oil varies greatly with temperature, and the water viscosity much less so. As for their densities, liquid water and oil are much less compressible and thermally expansive than other fluids such as gases. To create a higher pressure-temperature cross-coupling, we increase the compressibility coefficient $c$ and the thermal expansion coefficient $\beta$ of the oil density as defined in Section 2.1.4 of \cite{roy2019block}.

\begin{table}[!htb]\centering
\caption{Increased cross-coupling with $\Delta t = 2$ days. 2D Isotropic Well case. $N=160$. Average linear iterations per nonlinear iteration.}\label{tab:crosscoup1}
\begin{tabular}{|c c c c c c|}
 \hline
  Method/Factor & 1 & 5 & 10 & 15 & 20 \\
 \hline
   \rtwo{CPR-AMG}           & 8.29 & 8.29 & 8.42   & 8.63  & 8.8\\
   \rtwo{CPTR-Block(AMG)}   & 5.29 & 5.29 & 5.57   & 6	& 6.27\\
   \rtwo{CPTR-uAMG}    & 5.57 & 6.14 & 7	& 8.5	& 10.6\\
   \rtwo{CPTR-BD(LU)}    & 8.29 & 9.14 & 10.1   & 11.1  & 12.8  \\
   \rtwo{CPTR-BD(AMG)}   & 8.29 & 9.14 & 10.3	& 11.1	& 12.8\\
 \hline
 \end{tabular}
\end{table}

\begin{table}[!htb]\centering
\caption{Increased cross-coupling with $\Delta t = 0.1$ days. 2D Isotropic Well case. $N=160$. Average linear iterations per nonlinear iteration.}\label{tab:crosscoup2}
\begin{tabular}{|c c c c c c|}
 \hline
  Method/Factor & 1 & 15 & 30 & 35 & 40 \\
 \hline
   \rtwo{CPR-AMG}            & 5.33 & 6.17 & 5    & 5.17 & 6.33 \\
   \rtwo{CPTR-Block(AMG)}    & 5	 & 5.83 & 4    & 4.33 &	5.17 \\
   \rtwo{CPTR-uAMG}     & 4.83 & 8.33 & 21.7 & 92.8 &	$>$200 \\
   \rtwo{CPTR-BD(LU)}     & 5.83 & 7.83 & 12.3 & 19.7 & $>$200 \\
   \rtwo{CPTR-BD(AMG)}    & 6    & 7.5  & 12.3 & 19.7 & $>$200 \\
 \hline
 \end{tabular}
\end{table} 

We consider the 2D isotropic well case from Section \ref{sec:refine} with $N=160$. Both coefficients $c$ and $\beta$ are multiplied by the same increasing factor. Since larger coefficients result in a stiffer problem, smaller time-steps are needed for the convergence of Newton's method. Results for a time-step of 2 days and 0.1 days are shown in Tables \ref{tab:crosscoup1} and \ref{tab:crosscoup2}, respectively.
For CPTR, we observe that the block preconditioner exhibits a very small increase in iterations, while the other methods have a steady increase in iterations. It is clear that the block preconditioner is much better at dealing with a strong cross-coupling than the other versions of CPTR.  For CPR, it shows similar robustness with respect to a strong cross-coupling. Note that by taking smaller time-steps, we are in a regime where CPR is competitive. Since CPR performs appropriately but CPTR (other than \rtwo{CPTR-Block(AMG)}) does not, this may indicate that the coupling $A_{Tp}$ is more important than $A_{pT}$. The temperature solution provided in the first stage of CPTR by solvers other than the block preconditioner is not accurate. It remains to be seen if $A_{pT}$ is more important in other cases.

In Table \ref{tab:crosscoup2}, the number of iterations eventually blows up for uAMG, \rtwo{BD(LU)}, and \rtwo{BD(AMG)}. However, in the case of uAMG, this may not necessarily be because of a strong cross-coupling but perhaps that the properties of $A_{TT}$ are no longer amenable to the application of AMG. Indeed, the derivative of density with respect to temperature is negative. Consider the derivative of the oil enthalpy
\begin{equation}\label{eq:deriventhalpy}
 (c_{v,o} \rho_o  T )_T = c_{v,o}\rho_o + c_{v,o} (\rho_o)_T T,
\end{equation}
which is positive for the default parameters of the earlier sections. However, for the larger factors in Table \ref{tab:crosscoup2}, the negative influence of $(\rho_o) _T$ begins to dominate. For the largest factors, this results in $A_{TT}$ losing diagonal dominance, becoming indefinite, and eventually having a negative diagonal. This is, of course, far from the ideal properties for the convergence of AMG. For uAMG, the coarse pressure problem is coupled through temperature using the prolongation operator for the temperature coarse grid. Since that operator is not descriptive, this may indicate why uAMG is doing worse than \rtwo{BD(AMG)}. 

On the other hand, the block preconditioner does not apply AMG to $A_{TT}$, but rather to the approximate Schur complement $\tilde S_T$. Since this approximation does not include $(\rho_o) _T$, it will not suffer from a larger thermal expansivity coefficient. 


\section{Conclusion}\label{sec:conclusion}

In this work, we have implemented a fully implicit parallel non-isothermal multiphase flow in porous media simulator including two preconditioning strategies, CPR and CPTR. The first stage of CPTR can be a block preconditioner with our own Schur complement approximation, or an unknown-based AMG method. Standard decoupling operators are available for CPR, as well as their extension for CPTR. 

On coarse grids, both CPR and CPTR exhibit a similar number of iterations, which means that CPR outperforms CPTR in terms of computational cost. In these cases, heat diffusion is not very significant, so CPTR is not necessary. In other cases, however, CPTR displays much better scalability in terms of mesh refinement, as well as parallelization. Additionally, the first stage of CPTR is vulnerable to decoupling errors with saturation variables which can negate the scalability advantages of the method. Algebraic decoupling operators do not solve this issue. To deal with a strong pressure-temperature cross-coupling, our block preconditioner is better than the less coupled alternatives in the first stage of CPTR.

We believe that more work is needed to understand the decoupling error from the first stage of CPTR. To regain the scalability that CPTR displays in other cases, an accurate energy equation is key. It is still unclear if this can be done with another CPR-like multi-stage preconditioner.

%

\bibliographystyle{siamplain}
{
\interlinepenalty=10000
\bibliography{bibmultithermal}

\begin{thebibliography}{10}

\bibitem{alnaes2014unified}
{\sc M.~S. Aln{\ae}s, A.~Logg, K.~B. {\O}lgaard, M.~E. Rognes, and G.~N.
  Wells}, {\em Unified form language: A domain-specific language for weak
  formulations of partial differential equations}, ACM Transactions on
  Mathematical Software (TOMS), 40 (2014), p.~9.

\bibitem{petsc-web-page}
{\sc S.~Balay, S.~Abhyankar, M.~F. Adams, J.~Brown, P.~Brune, K.~Buschelman,
  L.~Dalcin, V.~Eijkhout, W.~D. Gropp, D.~Kaushik, M.~G. Knepley, L.~C.
  McInnes, K.~Rupp, B.~F. Smith, S.~Zampini, H.~Zhang, and H.~Zhang}, {\em
  {PETS}c {W}eb page}, 2017, \url{http://www.mcs.anl.gov/petsc}.

\bibitem{bennison1998prediction}
{\sc T.~Bennison}, {\em Prediction of heavy oil viscosity}, in Presented at the
  {IBC} Heavy Oil Field Development Conference, vol.~2, 1998, p.~4.

\bibitem{bui2017algebraic}
{\sc Q.~M. Bui, H.~C. Elman, and J.~D. Moulton}, {\em Algebraic multigrid
  preconditioners for multiphase flow in porous media}, SIAM Journal on
  Scientific Computing, 39 (2017), pp.~S662--S680.

\bibitem{bui2018algebraic}
{\sc Q.~M. Bui, L.~Wang, and D.~Osei-Kuffuor}, {\em Algebraic multigrid
  preconditioners for two-phase flow in porous media with phase transitions},
  Advances in water resources, 114 (2018), pp.~19--28.

\bibitem{chen2009well}
{\sc Z.~Chen and Y.~Zhang}, {\em Well flow models for various numerical
  methods}, International Journal of Numerical Analysis \& Modeling, 6 (2009).

\bibitem{christie2001tenth}
{\sc M.~A. Christie and M.~J. Blunt}, {\em Tenth {SPE} comparative solution
  project: A comparison of upscaling techniques}, in SPE Reservoir Simulation
  Symposium, Society of Petroleum Engineers, 2001.

\bibitem{clees2005amg}
{\sc T.~Clees}, {\em {AMG} strategies for {PDE} systems with applications in
  industrial semiconductor simulation}, PhD thesis, Universit\"{a}t zu
  K\"{o}ln, 2005.

\bibitem{clees2010efficient}
{\sc T.~Clees and L.~Ganzer}, {\em An efficient algebraic multigrid solver
  strategy for adaptive implicit methods in oil-reservoir simulation}, SPE
  Journal, 15 (2010), pp.~670--681.

\bibitem{coats1980situ}
{\sc K.~H. Coats}, {\em In-situ combustion model}, Society of Petroleum
  Engineers Journal, 20 (1980), pp.~533--554.

\bibitem{cremon2019multi}
{\sc M.~A. Cremon, N.~Castelletto, and J.~A. White}, {\em Multi-stage
  preconditioners for thermal-compositional-reactive flow in porous media},
  arXiv preprint arXiv:1912.04385,  (2019).

\bibitem{darcy1856fontaines}
{\sc H.~Darcy}, {\em Les fontaines publiques de la ville de {D}ijon}, Victor
  Dalmont, 1856.

\bibitem{elman2014finite}
{\sc H.~C. Elman, D.~J. Silvester, and A.~J. Wathen}, {\em Finite elements and
  fast iterative solvers: with applications in incompressible fluid dynamics},
  Oxford University Press, USA, 2014.

\bibitem{eymard2000finite}
{\sc R.~Eymard, T.~Gallou{\"e}t, and R.~Herbin}, {\em Finite volume methods},
  Handbook of Numerical Analysis, 7 (2000), pp.~713--1018.

\bibitem{falgout2002hypre}
{\sc R.~D. Falgout and U.~M. Yang}, {\em hypre: A library of high performance
  preconditioners}, in International Conference on Computational Science,
  Springer, 2002, pp.~632--641.

\bibitem{ml-guide}
{\sc M.~W. Gee, C.~M. Siefert, J.~J. Hu, R.~S. Tuminaro, and M.~G. Sala}, {\em
  {ML} 5.0 smoothed aggregation user’s guide}, Tech. Report SAND2006-2649,
  Sandia National Laboratories, 2006.

\bibitem{grabowski1979fully}
{\sc J.~W. Grabowski, P.~K. Vinsome, R.~C. Lin, G.~Behie, and B.~Rubin}, {\em A
  fully implicit general purpose finite-difference thermal model for in situ
  combustion and steam}, in SPE Annual Technical Conference and Exhibition,
  Society of Petroleum Engineers, 1979.

\bibitem{gries2015system}
{\sc S.~Gries}, {\em System-{AMG} Approaches for Industrial Fully and Adaptive
  Implicit Oil Reservoir Simulations}, PhD thesis, Universit{\"a}t zu K{\"o}ln,
  2015.

\bibitem{gries2014preconditioning}
{\sc S.~Gries, K.~St{\"u}ben, G.~L. Brown, D.~Chen, and D.~A. Collins}, {\em
  Preconditioning for efficiently applying algebraic multigrid in fully
  implicit reservoir simulations}, SPE Journal, 19 (2014), pp.~726--736.

\bibitem{henson2002boomeramg}
{\sc V.~E. Henson and U.~M. Yang}, {\em Boomer{AMG}: A parallel algebraic
  multigrid solver and preconditioner}, Applied Numerical Mathematics, 41
  (2002), pp.~155--177.

\bibitem{jenny2006adaptive}
{\sc P.~Jenny, S.~H. Lee, and H.~A. Tchelepi}, {\em Adaptive fully implicit
  multi-scale finite-volume method for multi-phase flow and transport in
  heterogeneous porous media}, Journal of Computational Physics, 217 (2006),
  pp.~627--641.

\bibitem{kell1975density}
{\sc G.~S. Kell}, {\em Density, thermal expansivity, and compressibility of
  liquid water from 0. deg. to 150. deg.. correlations and tables for
  atmospheric pressure and saturation reviewed and expressed on 1968
  temperature scale}, Journal of Chemical and Engineering Data, 20 (1975),
  pp.~97--105.

\bibitem{kirby2018solver}
{\sc R.~C. Kirby and L.~Mitchell}, {\em Solver composition across the
  {PDE}/linear algebra barrier}, SIAM Journal on Scientific Computing, 40
  (2018), pp.~C76--C98.

\bibitem{lacroix2003iterative}
{\sc S.~Lacroix, Y.~V. Vassilevski, J.~Wheeler, and M.~F. Wheeler}, {\em
  Iterative solution methods for modeling multiphase flow in porous media fully
  implicitly}, SIAM Journal on Scientific Computing, 25 (2003), pp.~905--926.

\bibitem{lacroix2000iterative}
{\sc S.~Lacroix, Y.~V. Vassilevski, and M.~F. Wheeler}, {\em Iterative solvers
  of the implicit parallel accurate reservoir simulator ({IPARS}), {I}: single
  processor case}, {TICAM} report 00-28, The University of Texas at Austin,
  2000.

\bibitem{leveque2002finite}
{\sc R.~J. LeVeque}, {\em Finite volume methods for hyperbolic problems},
  vol.~31, Cambridge University Press, 2002.

\bibitem{meijerink1977iterative}
{\sc J.~A. Meijerink and H.~A. van~der Vorst}, {\em An iterative solution
  method for linear systems of which the coefficient matrix is a symmetric
  {M}-matrix}, Mathematics of Computation, 31 (1977), pp.~148--162.

\bibitem{muskat1936flow}
{\sc M.~Muskat and M.~W. Meres}, {\em The flow of heterogeneous fluids through
  porous media}, Physics, 7 (1936), pp.~346--363.

\bibitem{patankar1980numerical}
{\sc S.~V. Patankar}, {\em Numerical heat transfer and fluid flow}, Series on
  Computational Methods in Mechanics and Thermal Science, Hemisphere Publishing
  Corporation (CRC Press, Taylor \& Francis Group), 1980.

\bibitem{peaceman1978interpretation}
{\sc D.~W. Peaceman}, {\em Interpretation of well-block pressures in numerical
  reservoir simulation (includes associated paper 6988)}, Society of Petroleum
  Engineers Journal, 18 (1978), pp.~183--194.

\bibitem{rathgeber2016firedrake}
{\sc F.~Rathgeber, D.~A. Ham, L.~Mitchell, M.~Lange, F.~Luporini, A.~T.~T.
  McRae, G.-T. Bercea, G.~R. Markall, and P.~H.~J. Kelly}, {\em Firedrake:
  automating the finite element method by composing abstractions}, ACM
  Transactions on Mathematical Software (TOMS), 43 (2016), p.~24.

\bibitem{riviere2008discontinuous}
{\sc B.~Rivi\`ere}, {\em Discontinuous {G}alerkin methods for solving elliptic
  and parabolic equations: theory and implementation}, SIAM, 2008.

\bibitem{roy2019preconditioning}
{\sc T.~Roy}, {\em Preconditioning for thermal reservoir simulation}, {DPhil}
  thesis, University of Oxford, 2019.

\bibitem{roy2019block}
{\sc T.~Roy, T.~B. J{\"o}nsth{\"o}vel, C.~Lemon, and A.~J. Wathen}, {\em A
  block preconditioner for non-isothermal flow in porous media}, Journal of
  Computational Physics, 395 (2019), pp.~636--652.

\bibitem{ruge1986amg}
{\sc J.~W. Ruge}, {\em {AMG} for problems of elasticity}, Applied Mathematics
  and Computation, 19 (1986), pp.~293--309.

\bibitem{ruge1987algebraic}
{\sc J.~W. Ruge and K.~St{\"u}ben}, {\em Algebraic multigrid}, in Multigrid
  methods, vol.~3 of Frontiers in Applied Mathematics, SIAM, Philadelphia,
  1987, ch.~4, pp.~73--130.

\bibitem{saad2003iterative}
{\sc Y.~Saad}, {\em Iterative methods for sparse linear systems}, SIAM, 2003.

\bibitem{saad1986gmres}
{\sc Y.~Saad and M.~H. Schultz}, {\em {GMRES}: A generalized minimal residual
  algorithm for solving nonsymmetric linear systems}, SIAM Journal on
  Scientific and Statistical Computing, 7 (1986), pp.~856--869.

\bibitem{sahni2000electromagnetic}
{\sc A.~Sahni, M.~Kumar, and R.~B. Knapp}, {\em Electromagnetic heating methods
  for heavy oil reservoirs}, in SPE/AAPG Western Regional Meeting, Society of
  Petroleum Engineers, 2000.

\bibitem{scheichl2003decoupling}
{\sc R.~Scheichl, R.~Masson, and J.~Wendebourg}, {\em Decoupling and block
  preconditioning for sedimentary basin simulations}, Computational
  Geosciences, 7 (2003), pp.~295--318.

\bibitem{stuben2001introduction}
{\sc K.~St{\"u}ben}, {\em An introduction to algebraic multigrid}, in
  Multigrid, U.~Trottenberg, C.~W. Oosterlee, and A.~Schuller, eds., Academic
  Press, 2001, pp.~413--532.

\bibitem{wallis1983incomplete}
{\sc J.~R. Wallis}, {\em Incomplete {G}aussian elimination as a preconditioning
  for generalized conjugate gradient acceleration}, in SPE Reservoir Simulation
  Symposium, Society of Petroleum Engineers, 1983.

\bibitem{wallis1985constrained}
{\sc J.~R. Wallis, R.~P. Kendall, and T.~E. Little}, {\em Constrained residual
  acceleration of conjugate residual methods}, in SPE Reservoir Simulation
  Symposium, Society of Petroleum Engineers, 1985.

\bibitem{wang2017multigrid}
{\sc L.~Wang, D.~Osei-Kuffuor, R.~Falgout, I.~Mishev, and J.~Li}, {\em
  Multigrid reduction for coupled flow problems with application to reservoir
  simulation}, in SPE Reservoir Simulation Conference, Society of Petroleum
  Engineers, 2017.

\bibitem{wong2018sequential}
{\sc Z.~Y. Wong}, {\em Sequential-implicit Newton's Method for Geothermal
  Reservoir Simulation}, PhD thesis, Stanford University, 2018.

\end{thebibliography}
}

\appendix

\section{Physical quantities}\label{sec:murho}
For the oil viscosity $\mu_o$, we choose a correlation from \cite{bennison1998prediction} where $\mu_o$ varies hugely with temperature. For the oil density $\rho_o$, we choose a function of pressure and temperature given by an exponential correlation. We choose an ``API gravity'' of 10, such that the density of oil and water are the same at surface conditions. Details and choices of parameters for $\mu_o$ and $\rho_o$ can be found in Section 2.1.4 of \cite{roy2019block}.

For the water viscosity $\mu_w$, we use the following correlation \cite{grabowski1979fully}:
\begin{equation}\label{eq:watermu}
 \mu_w (T_\mathrm{F}) = \frac{A}{-1 + B T_\mathrm{F} + C T_\mathrm{F}^2},
\end{equation}
which takes temperature $T_\mathrm F$ in $^\circ$F and returns viscosity in cp (0.001~kg\,m$^{-1}$\,s$^{-1}$). The parameters $A,B,C$ can be found in Table \ref{tab:watermu}. 

\begin{table}[hbt]
\centering
\caption{Parameters for the water viscosity correlation}\label{tab:watermu}
 \begin{tabular}{c c c}
 \hline
  $A$ & $B$ & $C$  \\
  \hline
  2.1850 &  0.04012 & 5.1547$\times10^{-6}$\\
  \hline
 \end{tabular}
\end{table}

For the water density, we use Trangenstein's modification of Kell's formula \cite{kell1975density} 
\begin{equation}
\rho_w (p,T_\mathrm{C}) =  \frac{E_0 + E_1 T + E_2 T^2 + E_3 T^3 + E_4 T^4 + E_5 T^5}{1 + E_6 T} e^{C_w (p-E_7)},
\end{equation}
which takes temperature $T_\mathrm{C}$ in $^\circ$C, and pressure $p$ in MPa. The parameters are given in Table \ref{tab:waterrho}.

\begin{table}[hbt!]
\centering
\caption{Parameters for the water density correlation}\label{tab:waterrho}
  \begin{tabular}{l l }
 \hline
  $E_0$ & 999.83952 \\
  $E_1$ & 16.955176 \\
  $E_2$ & -7.987$\times 10^{-3}$ \\
  \hline
 \end{tabular}
  \begin{tabular}{l l}
 \hline
 $E_3$ & -46.170461$\times 10^{-6}$ \\
 $E_4$ & 105.56302$\times 10^{-9}$ \\
 $E_5$ &-280.54353$\times 10^{-12}$ \\
  \hline
 \end{tabular}
   \begin{tabular}{l l}
 \hline
 $E_6$ & 16.87985$\times 10^{-3}$ \\
 $E_7$ & 10.2 \\
 $C_w$ & 3.98854$\times 10^{-4}$ \\
  \hline
 \end{tabular}
\end{table}

For the relative permeability, we simply set it the be equal to the phase saturation. If we set $S_w = 1 - S_o$, we get $k_{ro} = S_o$, and $k_{rw} = 1 - S_o$. 
\end{document}